\documentclass[12pt]{amsart}
\usepackage{fullpage}
\usepackage{enumitem}
\usepackage{tikz}
\usetikzlibrary{arrows,calc,positioning,through,intersections,backgrounds,matrix,fit,shapes,decorations}

\pgfdeclarelayer{backbackground}
\pgfdeclarelayer{background}
\pgfdeclarelayer{foreground}
\pgfsetlayers{backbackground,background,main,foreground}		
\usepackage{xcolor}
\usepackage{ifthen}
\usepackage{hhline}
\usepackage{caption}

\tikzset{
  vertex/.style={circle,thick,color=white,fill=black,minimum size=6pt, inner sep=3pt},
  edge/.style={line width=1.5pt, black},
  bigedge/.style={line width=3pt, black},
}
\newtheorem{theorem}{Theorem}[section]
\newtheorem{lemma}[theorem]{Lemma}

\newtheorem{observation}[theorem]{Observation}
\newtheorem{fact}[theorem]{Fact}
\newtheorem{example}[theorem]{Example}
\newtheorem{conjecture}[theorem]{Conjecture}
\newtheorem{claim}[theorem]{Claim}
\newtheorem*{claim*}{Claim}

\theoremstyle{definition}
\newtheorem{definition}[theorem]{Definition}
\newtheorem{remark}[theorem]{Remark}

\newcommand{\cA}{\mathcal{A}}

\newcommand{\cC}{\mathcal{C}}

\newcommand{\cG}{\mathcal{G}}

\newcommand{\cK}{\mathcal{K}}

\newcommand{\cP}{\mathcal{P}}

\newcommand{\cV}{\mathcal{V}}
\newcommand{\cW}{\mathcal{W}}
\newcommand{\cX}{\mathcal{X}}

\newcommand{\cZ}{\mathcal{Z}}

\newcommand{\floor}[1]{\left\lfloor #1 \right\rfloor}
\newcommand{\ceiling}[1]{\left\lceil #1 \right\rceil}
\newcommand{\eps}{\varepsilon}

\newcommand{\E}{\mathbb{E}}

\title{Powers of Hamiltonian cycles in multipartite graphs}
\author{Louis DeBiasio$^{1}$, Ryan Martin$^{2}$, Theodore Molla$^{3}$}

\begin{document}

\begin{abstract}
We prove that if $G$ is a $k$-partite graph on $n$ vertices in which all of the parts have order at most $n/r$ and every vertex is adjacent to at least a $1-1/r+o(1)$ proportion of the vertices in every other part, then $G$ contains the $(r-1)$-st power of a Hamiltonian cycle.   
\end{abstract}

\maketitle
\noindent\footnotetext[1]{Department of Mathematics, Miami University {\tt debiasld@miamioh.edu}. Research supported in part by Simons Foundation Collaboration Grant \# 283194 and NSF grant DMS-1954170.}

\noindent\footnotetext[2]{Department of Mathematics, Iowa State University {\tt rymartin@iastate.edu}. Research supported in part by Simons Foundation Collaboration Grants \# 353292 and \#709641.}

\noindent\footnotetext[3]{Department of Mathematics and Statistics, University of South Florida {\tt molla@usf.edu}. 
Research supported in part by 
NSF Grants DMS-1500121
and DMS-1800761.}

\section{Introduction}

For graphs $G$ and $H$, we say that $G$ has a perfect $H$-tiling if $G$ contains $|V(G)|/|V(H)|$ vertex disjoint copies of $H$.  For a positive integer $r$, the $r$-th power of $H$ denoted $H^r$, is the graph on $V(H)$ where $uv\in E(H^r)$ if and only if the distance between $u$ and $v$ in $H$ is at most $r$.  We refer to the $(r-1)$-st power of a cycle as an $(r-1)$-cycle.

Hajnal and Szemer\'edi \cite{HS} proved that for all positive integers $r$ and $n$, if $r$ divides $n$ and $G$ is a graph on $n$ vertices with $\delta(G)\geq \left(1-\frac{1}{r}\right)n$, then $G$ contains a perfect $K_r$-tiling.  Koml\'os, S\'ark\"ozy, and Szemer\'edi \cite{KSS4} proved that for all $r\geq 2$, there exists $n_0$ such that if $G$ is a graph on $n\geq n_0$ vertices with $\delta(G)\geq \left(1-\frac{1}{r}\right)n$, then $G$ contains a Hamiltonian $(r-1)$-cycle.  Note that if $r$ divides $n$ and $G$ contains a Hamiltonian $(r-1)$-cycle, then $G$ contains a perfect $K_r$-tiling, so the result of Koml\'os, S\'ark\"ozy, and Szemer\'edi is stronger for fixed $r$ and large $n$.

  A graph $G$ is a \textit{$k$-partite graph with ordered partition $\cP = (V_1, \dotsc, V_k)$}, if 
  $\cP$ is a partition of $V(G)$ and $V_i$ is an independent set for every $i \in [k]$. 
  For all $i\neq j\in [k]$, let 
  $$\delta_{ij}(G)=\frac{\min\{\deg_G(v, V_j): v\in V_i\}}{|V_j|} \qquad\text{ and }\qquad \delta_{\cP}(G)=\min_{i\neq j\in [k]} \delta_{ij}(G).$$ 

Fisher \cite{F} conjectured an analogue of the Hajnal-Szemer\'edi theorem in balanced multipartite graphs; that is, if $G$ is a balanced $r$-partite graph on $n$ vertices with \[
  \delta_{\cP}(G) \geq 1 - \frac{1}{r},
\]
then $G$ contains a perfect $K_r$-tiling.  An earlier example of Catlin \cite{C} provides a counterexample to Fisher's conjecture when $r$ is odd, but Magyar and Martin \cite{MM} proved that for $r=3$, Catlin's counterexample is the only one.  Then Martin and Szemer\'edi \cite{MS} proved Fisher's conjecture for $r=4$.  After a relatively large gap in activity, Keevash and Mycroft \cite{KM1} and independently Lo and Markstr\"om \cite{LM2} proved that for all $\gamma>0$ and $r\geq 2$, there exists $n_0$ such that for all $n\geq n_0$ in which $r$ divides $n$, if   $G$ is a balanced $r$-partite graph on $n$ vertices with \[
  \delta_{\cP}(G) \geq 1 - \frac{1}{r} + \gamma,
\]
then $G$ contains a perfect $K_r$-tiling.  Later, an exact version was proved by Keevash and Mycroft \cite{KM2} which again shows that Fisher's conjecture holds for sufficiently large $n$ unless $r$ is odd in which case Catlin's counterexample is the only one.

Our main result can be viewed as a strengthening of the asymptotic versions of all of the above results (both in the multipartite setting and in the ordinary setting).

\begin{theorem}\label{thm:main}
For all $k\ge r\ge 2$ and all $0<\gamma\leq \frac{1}{r}$, there exists $n_0$ such that for all $n \ge n_0$
the following holds. 
If $G$ is a $k$-partite graph on $n$ vertices with 
ordered partition $\cP = (V_1, \dotsc, V_k)$
such that $|V_i| \le n/r$ for all $i \in [k]$ and
\[
  \delta_{\cP}(G) \geq 1 - \frac{1}{r} + \gamma,
\]
then $G$ contains a Hamiltonian $(r-1)$-cycle.  
\end{theorem}

Note that the condition $|V_i|\leq n/r$ for all $i\in [k]$ is necessary for the existence of a Hamiltonian $(r-1)$-cycle since the $(r-1)$-st power of a cycle on $n$ vertices has independence number $\floor{n/r}$.  Also this result is seen to be asymptotically best possible by taking a complete $k$-partite graph with ordered partition $\cP=(V_1, \dots, V_k)$ and letting $V_i'\subseteq V_i$ for all $i\in [k]$ with $|V_i'|=\floor{|V_i|/r}+1$ and deleting all edges inside $V_1'\cup \dots \cup V_k'$ to get a $k$-partite graph $G$ with $\delta_{\cP}(G)$ just below $1-\frac{1}{r}$ which has independence number larger than $\floor{n/r}$ and thus does not contain a Hamiltonian $(r-1)$-cycle.

\section{Observations, definitions and tools}

\begin{observation}\label{obs:combine}
  It suffices to prove Theorem~\ref{thm:main} in the cases where
  $r\leq k\leq 2r-1$ and all of the parts have order at least $\frac{\gamma}{2r} n$. 
\end{observation}

\begin{proof}
  Suppose Theorem~\ref{thm:main} is true provided $2\leq r\leq k\leq 2r-1$ and $|V_i|\geq \frac{\gamma}{2r} n$ for all $i\in [k]$.  Now suppose for contradiction that there exists a counterexample to Theorem \ref{thm:main}.  Let $k'$ be minimal such that a counterexample exists.  Let $n_0$ be the value coming from Theorem \ref{thm:main} when $k=k'-1$ and $\gamma'=\frac{\gamma}{2r}$.  Let $G'$ be a $k'$-partite counterexample on $n\geq n_0$ vertices with ordered partition $\cP'=(U_1, \dots, U_{k'})$ where $k'$ is minimal.  
  
  We first claim that for all distinct $i, j\in [k']$, $|U_i|+|U_j|>n/r$.  Suppose not and without loss of generality suppose that $i=k'-1$ and $j=k'$; that is, suppose $|U_{k'-1}|+|U_{k'}|\leq n/r$.  Let $V_i=U_i$ for all $i\in [k'-2]$ and $V_{k'-1}=U_{k'-1}\cup U_{k'}$ and let $G$ be the $(k'-1)$-partite graph with ordered partition $\cP=(V_1, \dots, V_{k'-1})$ obtained by deleting all edges between $U_{k'-1}$ and $U_{k'}$.  Since $\deg_{G'}(v, V_{k'-1})\geq (1-\frac{1}{r}+\gamma)|U_{k'-1}|+(1-\frac{1}{r}+\gamma)|U_{k'}|=(1-\frac{1}{r}+\gamma)|U_{k'-1}\cup U_{k'}|$ for all $v\in V(G)\setminus V_{k'-1}$ we have 
  \begin{equation*}
    \delta_{\cP}(G) \geq 1-\frac{1}{r}+\gamma.
  \end{equation*}
  But now by minimality, $G\subseteq G'$ has a Hamiltonian $(r-1)$-cycle contradicting the fact that $G'$ does not.  
  Thus we may assume that $r\leq k'\leq 2r-1$ as otherwise the two smallest parts add up to at most $n/r$.  
  
  Now suppose $G$ has a part of order less than $\gamma'n=\frac{\gamma }{2r}n$; without loss of generality, suppose it is $U_{k'}$.  Because $|U_i| \le n/r$ for every $i \in [k']$ and $\sum_{i \in [k']} |U_i| = n$, the fact that $|U_{k'}| \le \gamma' n< n/r$ implies that $k'>r$.  By the above, we may suppose that all other parts have order greater than $\frac{n}{r}-\gamma' n$.  Now partition $U_{k'}$ arbitrarily as $\{U_1', \dots, U_{k'-1}'\}$ (allowing for empty sets in the partition) subject to $|U_i|+|U_i'|\leq n/r$ for all $i\in [k'-1]$.  Let $G$ be the $(k'-1)$-partite graph with ordered partition $\cP=(V_1, \dots, V_{k'-1})$ where $V_i=U_i\cup U_i'$ for all $i\in [k'-1]$.  Since $$\left(1-\frac{1}{r}+\gamma\right)|U_i|\geq \left(1-\frac{1}{r}+\gamma'\right)(|U_i|+\gamma' n) \geq \left(1-\frac{1}{r}+\gamma'\right)|V_i|,$$ we have 
   \begin{equation*}
    \delta_{\cP}(G) \geq 1-\frac{1}{r}+\gamma',
  \end{equation*}
  and thus by minimality and the choice of $n_0$, $G\subseteq G'$ has a Hamiltonian $(r-1)$-cycle contradicting the fact that $G'$ does not.
\end{proof}

The following simple fact is used implicitly throughout the paper.
\begin{fact}\label{fact:almost_spanning_degree}
  Let $\sigma > 0$ and $G$ be a $k$-partite graph on
  $n$ vertices with ordered partition 
  $\cP = (V_1, \dotsc, V_k)$ such that
  every part has order at least $\sigma n$.
  For every $U \subseteq V(G)$ such that $|U| \le \sigma^2 n$, 
  if $G' = G - U$, then $\delta_{\cP}(G') \ge \delta_{\cP}(G) - \sigma$.
\end{fact}
\begin{proof}
  For distinct $i,j \in [k]$ and every $v \in V(G') \cap V_i$, we have 
  \begin{equation*}
    \frac{\deg_{G'}(v, V(G') \cap V_j)}{|V(G') \cap V_j|} \ge 
    \frac{\deg_{G'}(v, V(G') \cap V_j)}{|V_j|} \ge 
    \frac{\deg_G(v, V_j)}{|V_j|} - \frac{|U|}{|V_j|}
    \ge \delta_{\cP}(G) - \sigma. \qedhere
  \end{equation*}
\end{proof}

\begin{definition}[$(r-1)$-path/$(r-1)$-walk]
  Let $G$ be a graph
  and let $\cW = x_1, \dotsc, x_\ell$ be an ordered
  sequence of vertices of $G$.
  The sequence $\cW$ is an \textit{$(r-1)$-walk of length $\ell$}
  if every $r$ consecutive vertices in $\cW$ form a clique in $G$. 
  If $\cW$ is an $(r-1)$-walk of length $\ell$,
  then it is an \textit{$(r-1)$-path of length $\ell$} 
  if there are no repeated vertices in the sequence $x_1, \dotsc, x_\ell$. 
\end{definition}

The following fact is immediate when one first observes that 
the number of $(r-1)$-walks of length $\ell$ that are not $(r-1)$-paths is at most
$\binom{\ell}{2} \cdot n^{\ell - 1}$,
and that, for every set $U \subseteq V(G)$, 
the total number of $(r-1)$-walks of length $\ell$ that contain a vertex from $U$ is 
at most $\ell \cdot |U| \cdot n^{\ell - 1}$.
Throughout the remainder of the proof, 
we use the notation $a \ll b$ to indicate that there exists an increasing
function $f(b)$ such that the result holds for every $a \le f(b)$.

\begin{fact}\label{fact:walks_to_paths}
  Suppose $\frac{1}{n} \ll \sigma \ll \alpha, \frac{1}{\ell}$ and let $G$ be an $n$-vertex graph
  and $U \subseteq V(G)$ where $|U| \le \sigma n$.
  If $\cW$ is a collection of at least $(\alpha n)^\ell$ $(r-1)$-walks of length $\ell$, 
  then at least $(\sigma n)^\ell$ of the walks in $\cW$ are $(r-1)$-paths 
  that avoid the set $U$.
\end{fact}

To motivate the following definition, let us first comment that, at various times, we will need to 
connect disjoint $(r-1)$-paths to form longer $(r-1)$-paths.
To highlight some issues that might arise in as simple a setting as possible, consider the case when $k=r=3$
and let $G$ be a balanced $3$-partite graph with ordered partition $(V_1, V_2, V_3)$ and let
$P_1 = u_1,\dotsc,u_6$ and $P_2=w_1,\dotsc,w_6$ be two disjoint $2$-paths each on $6$ vertices.
Suppose that we would like to find a $2$-path $Q$ so that the sequence $P_1QP_2$ 
is itself a $2$-path.
This would be impossible if, say, $u_4 \in V_1$, $u_5 \in V_2$, and $u_6 \in V_3$
while $w_1 \in V_2$, $w_2 \in V_1$ and $w_3 \in V_3$.
To see this, note that, in this setting, if 
$u_4 \in V_1$, $u_5 \in V_2$, and $u_6 \in V_3$ and $u_1, \dotsc, u_{3p}$ is a $2$-path,
then for every $0 \le i \le p - 1$ and $j \in [3]$, we must have that $u_{3i + j} \in V_j$.
To deal with issues such as this, we will require that $(r-1)$-walks conform to the following
definition.
\begin{definition}[Properly terminated]
Suppose that $G$ is a $k$-partite graph with ordered partition $(V_1, \dotsc, V_k)$ and
let $W=v_1v_2\dots v_p$ be an $(r-1)$-walk where $p \ge r$.
We say that $W$ is \textit{properly terminated} if $v_i\in V_i$ and $v_{p-r+i}\in V_i$ for all $i \in [r]$.
That is, $W$ is properly terminated if its first $r$ vertices traverse the sets $V_1, \dotsc, V_r$
in order and its last $r$ vertices traverse the sets $V_1, \dotsc, V_r$ in order.

More generally, if $\mathcal{P} = (U_1, \dotsc, U_r)$ is an ordered sequence of $r$ disjoint sets, we say that the \textit{initial $r$ vertices of $W$ respect the sequence $\mathcal{P}$} if $v_i \in U_i$ for every $i \in [r]$.
Similarly, we say that the \textit{final $r$ vertices of $W$ respect the sequence $\mathcal{P}$} if $v_{p - r +  i} \in U_i$ for every $i \in [r]$.
So, $W$ is properly terminated if both the initial $r$ vertices of $W$ and
the final $r$ vertices of $W$ respect the sequence $(V_1, \dotsc, V_r)$.
\end{definition}

\begin{definition}[Balanced]
  Let $\cP$ be a collection of disjoint sets.
  We say that $\cP$ is \textit{balanced} if every set in $\cP$ has the
  same order.
  
  If $G$ is an $r$-partite graph with ordered partition 
  $\cP = (V_1, \dotsc, V_r)$, we say that $G$ is 
  \textit{balanced} if $\cP$ is balanced
  and we say that a set $U\subseteq V(G)$ is 
  \textit{balanced} if $|U\cap V_i|=|U\cap V_j|$ for all $i,j\in [r]$.
\end{definition}

A few times in the proof we will 
make use of a Chernoff bound on the concentration of binomial and 
hypergeometric distributions \cite[Corollary 2.3 and Theorem 2.10]{JLR}
\begin{theorem}[Chernoff bound]\label{chernoff}
Suppose $X$ has binomial or hypergeometric distribution and $0<a<3/2$. Then
$\mathbb{P}(|X - \E X| \ge a\E X) \le 2
e^{-\frac{a^2}{3}\E X}$.
\qed
\end{theorem}

\section{Overview of the proof}

We are attempting to prove that all sufficiently large $k$-partite graphs, in which all parts have at most $n/r$ vertices, with proportional minimum degree at least $1-\frac{1}{r}+\gamma$ have a Hamiltonian $(r-1)$-cycle.  We are able to split the work into two tasks. 

The first (and main) task is to prove the result in the case of \emph{balanced} $r$-partite graphs. Lemma~\ref{lem:balanced} below establishes that in a large balanced $r$-partite graph, and two properly terminated $(r-1)$-paths with the same ordering, $K$ and $K'$, there is a Hamiltonian $(r-1)$-path that starts with $K$ and ends with $K'$. If the graph is balanced and $r$-partite, then we simply apply this with $K=K'$ and we are done. If not, then we use Lemma~\ref{lem:seq} below to partition the graph into balanced $r$-partite pieces and then stitch them together to create the $(r-1)$-cycle we require. 

\begin{lemma}[Balanced case]\label{lem:balanced}
  For every $r \ge 2$ and $\gamma < \frac{1}{r}$, there exists $n_0$
  such that for every $n \ge n_0$ the following holds.
  Let $G$ be a balanced $r$-partite graph 
  on $n$ vertices
  with ordered partition $\cP = (V_1, \dotsc, V_r)$
  such that
  \begin{equation*}
    \delta_{\cP}(G) \ge 1 - \frac{1}{r} + \gamma.
  \end{equation*}
  Suppose that $K$ and $K'$ are $r$-cliques such that either $K = K'$ or $K \cap K' = \emptyset$
  and let $v_i := V_i \cap K$ and $v'_i := V_i \cap K'$ for every $i \in [r]$.
  Then there is a Hamiltonian $(r-1)$-path $P$ of $G - (K \cup K')$ such that
  $v_1, \dotsc, v_r,P,v'_1, \dotsc, v_r'$ is an $(r-1)$-walk in $G$. 
\end{lemma}

The second task is to show that $G$ can be partitioned into a small number of balanced $r$-partite graphs that each contain a Hamiltonian $(r-1)$-path and that these $(r-1)$-paths can be stitched together to form a Hamiltonian $(r-1)$-path of the original graph $G$. Lemma~\ref{lem:seq} below shows that the graph can be partitioned into balanced $r$-partite graphs $G_1, \dots, G_\ell$, each with the appropriate minimum degree condition, together with short $(r-1)$-paths connecting $G_i$ to $G_{i+1}$ in sequence in such a way that every vertex is accounted for.  Then applying Lemma~\ref{lem:balanced} to each $G_i$ we will construct the desired Hamiltonian $(r-1)$-cycle.

The technical issue for finding the partition is essentially numerical, requiring the sizes of the sets forming each $G_i$ to be the same and to partition each vertex class. Once these constraints are achieved, we are able to meet the minimum degree condition by applying a Chernoff bound to show that a randomly chosen partition satisfying the numerical constraints will have the required degree condition with high probability.

\begin{figure}
    \centering
    \def\xgap{2.7}
    \def\ygap{1.75}
    \def\wid{2.2}
    \begin{tikzpicture}
        \begin{pgfonlayer}{foreground}
            \foreach \j in {2,3,5,6}{
                \node[vertex] (1;\j) at (\j*\xgap,1*\ygap-0.7) {};
            }
            \foreach \j in {1,3,4,6}{
                \node[vertex] (-1;\j) at (\j*\xgap,-1*\ygap+0.7) {};
            }
        \end{pgfonlayer}
        \begin{pgfonlayer}{foreground}
            \node at (1*\xgap,2*\ygap) {$V_1$};
            \node at (2*\xgap,2*\ygap) {$V_2$};
            \node at (3*\xgap,2*\ygap) {$V_3$};
            \node at (4*\xgap,2*\ygap) {$V_4$};
            \node at (5*\xgap,2*\ygap) {$V_5$};
            \node at (6*\xgap,2*\ygap) {$V_6$};
            \node at (6*\xgap+1.8,\ygap) {$G_j$};
            \node at (6*\xgap+2.0,-\ygap) {$G_{j+1}$};
            \node at (2.5*\xgap,0.2) {$P_j$};
        \end{pgfonlayer}
        \begin{pgfonlayer}{main}
            \draw [edge] (1;2) -- (1;3);
            \draw [edge] (1;3) -- (1;5);
            \draw [edge] (1;5) -- (1;6);
            \draw [edge] (1;2) to[out=7, in=173] (1;5);
            \draw [edge] (1;3) to[out=7, in=173] (1;6);
            \draw [edge] (1;2) to[out=10, in=170] (1;6);
            \draw [edge] (-1;1) -- (-1;3);
            \draw [edge] (-1;3) -- (-1;4);
            \draw [edge] (-1;4) -- (-1;6);
            \draw [edge] (-1;1) to[out=-7, in=-173] (-1;4);
            \draw [edge] (-1;3) to[out=-7, in=-173] (-1;6);
            \draw [edge] (-1;1) to[out=-7, in=-173] (-1;6);
            \draw [edge] (1;3) -- (-1;1);
            \draw [edge] (1;5) -- (-1;1);
            \draw [edge] (1;5) -- (-1;3);
            \draw [edge] (1;6) -- (-1;1);
            \draw [edge] (1;6) -- (-1;3);
            \draw [edge] (1;6) -- (-1;4);
            \foreach \j in {2,3,5,6}{
                \node at ($(1;\j)+(0,1.1)$) {$V({\j, j})$};
            }
            \foreach \j in {1,3,4,6}{
                \node at ($(-1;\j)+(0,-1.1)$) {$V({\j, j+1})$};
            }
        \end{pgfonlayer}
        \begin{pgfonlayer}{background}
            \foreach \j in {2,3,5,6}{
                \draw [rounded corners, black, fill=white] (\j*\xgap-\wid/2,1*\ygap-1.0) rectangle (\j*\xgap+\wid/2,1*\ygap+0.8);
            }
            \foreach \j in {1,3,4,6}{
                \draw [rounded corners, black, fill=white] (\j*\xgap-\wid/2,-1*\ygap-0.8) rectangle (\j*\xgap+\wid/2,-1*\ygap+1.0);
            }
            \draw [dotted,rounded corners, black, very thick] (2*\xgap-1.4,0*\ygap+0.6) rectangle (6*\xgap+1.4,2*\ygap-0.8);
            \draw [dotted,rounded corners, black, very thick] (1*\xgap-1.4,0*\ygap-0.6) rectangle (6*\xgap+1.4,-2*\ygap+0.8);
        \end{pgfonlayer}
        \begin{pgfonlayer}{backbackground}
            \foreach \j in {1,2,3,4,5,6}{
                \draw [rounded corners, white, fill=gray!40] (\j*\xgap-\wid/2,-1.5*\ygap-0.6) rectangle (\j*\xgap+\wid/2,1.5*\ygap+0.6);
            }
        \end{pgfonlayer}
    \end{tikzpicture}
    \caption{An example for Lemma~\ref{lem:seq} in the case where $k=6$, $r=4$, $(i_{j,1},i_{j,2},i_{j,3},i_{j,4})=(2,3,5,6)$, $(i_{j+1,1},i_{j+1,2},i_{j+1,3},i_{j+1,4})=(1,3,4,6)$ and the $8$-vertex $3$-path, $P_j$.}
    \label{fig:seq}
\end{figure}
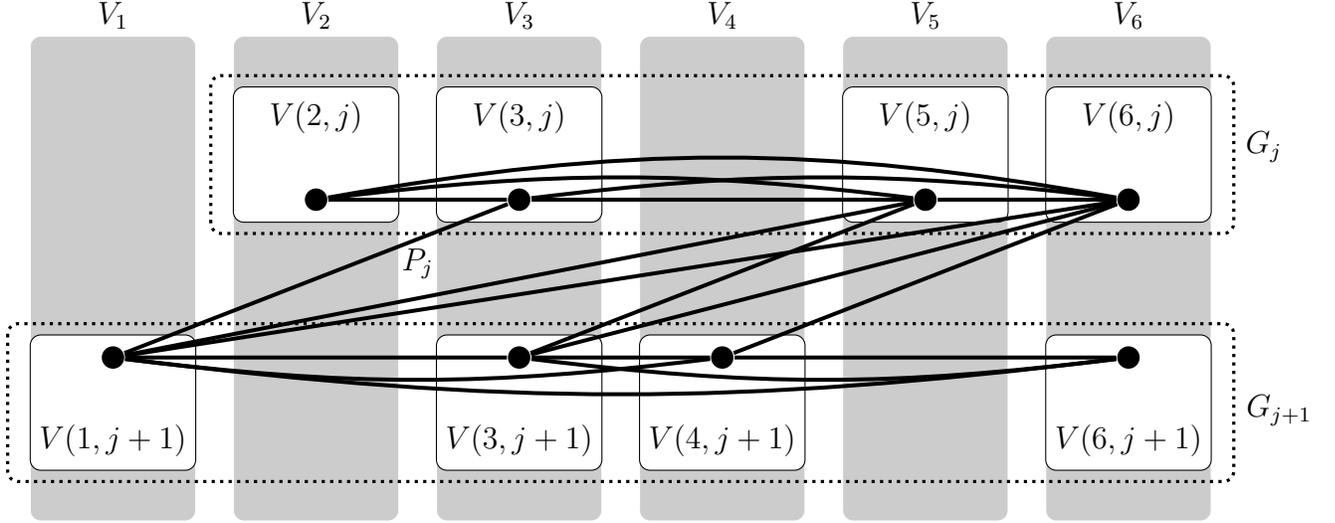

\begin{lemma}[Partitioning and Sequencing]\label{lem:seq}
For all $r\geq 2$, $0<\gamma\leq \frac{1}{r}$, and $r < k\leq 2r-1$, there exist constants $0<\frac{1}{n_0}\ll \beta\ll \sigma\ll \gamma$ such that
if $G$ is a $k$-partite graph on $n\geq n_0$ vertices with ordered partition 
$\cP = (V_1, \dotsc, V_k)$ in which $\gamma n\le |V_k| \le |V_{k-1}| \le \dotsm \le |V_1| \le \frac{n}{r}$
and 
$$
  \delta_{\cP}(G) \ge 1 - \frac{1}{r} + \gamma,
$$
then
there exists an $(r-1)$-path $P_0'$ with $|V(P_0')|\leq \beta n$ such that if $V_i'=V_i\setminus V(P_0')$ for $i \in [k]$, then the following holds:
\begin{enumerate}[label=(A\arabic*)]
  \item\label{A1} there exists a positive integer $\ell$ such that for all $i\in [k]$, there exists a partition $V_i'=\{V({i,1}), \dots, V({i,\ell})\}$ (with $V({i,j})$ possibly empty) such that for all $j\in [\ell]$ there exists $1\leq i_{j,1}<\dots< i_{j,r}\leq k$ such that $|V({i_{j,1}, j})|=\dots=|V({i_{j,r}, j})|\geq \beta n$ and if $i \in [k] \setminus \{i_{j,1}, \dotsc, i_{j,r} \}$, then
    $V(i,j) = \emptyset$, and
  \item\label{A2} letting $\cP_j=(V({i_{j,1},j}), \dots, V({i_{j,r},j}))$ and $G_j$ be the natural $r$-partite graph induced by $\cP_j$, we have that $\delta_{\cP_j}(G_j)\ge 1-\frac{1}{r}+\frac{\gamma}{2}$.
  \item\label{A3}
    We can prepend $r$ vertices and append $r$ vertices to $P'_0$ to create an
    $(r-1)$-path $P_0$ such that 
    the initial $r$ vertices of $P_0$ respect the sequence $\mathcal{P}_\ell$
    and the final $r$ vertices of $P_0$ respect the sequence $\mathcal{P}_1$.
  \item\label{A4} There exist vertex disjoint $(r-1)$-paths $P_1, \dots, P_{\ell-1}$ in $G - V(P_0)$ each on $2r$ vertices such that for all $j \in [\ell-1]$ the initial $r$ vertices of $P_j$ respect the sequence $\mathcal{P}_j$ and the final $r$ vertices of $P_j$ respect the sequence $\mathcal{P}_{j+1}$ (Fig.~\ref{fig:seq}).
\end{enumerate}
\end{lemma}

Lemma~\ref{lem:balanced} and Lemma~\ref{lem:seq} together immediately imply Theorem \ref{thm:main}.

\begin{proof}[Proof of Theorem \ref{thm:main}]
By Lemma~\ref{lem:balanced} we can assume
$k > r$ and 
by Observation~\ref{obs:combine} we can assume 
that $k\leq 2r-1$ and
every part of $\cP$ has order at least $\gamma' n$
where $\gamma \ge \gamma' > 0$.  Without loss of
generality we can further assume that
$\gamma' n \le |V_k| \le |V_{k-1}| \le \dotsm \le |V_1| \le \frac nk$.
Therefore, we can apply Lemma \ref{lem:seq} to $G$ (with $\gamma'$ playing the role of $\gamma$).  Define $P_\ell = P_0$ and for each $G_j$, apply Lemma~\ref{lem:balanced} to $G_j$ with $K=P_{j-1}\cap G_j$ and $K'=P_{j}\cap G_j$ to get a Hamiltonian $(r-1)$-path $Q_i$.  Now $P_0Q_1\dots Q_\ell$ is the desired Hamiltonian $(r-1)$-cycle.  
\end{proof}

In Section \ref{sec:lems} we describe the three lemmas needed to prove Lemma~\ref{lem:balanced}.  Then in Sections \ref{sec:con} to \ref{sec:cov}, we prove those lemmas.  Finally in Section \ref{sec:seq} we prove Lemma \ref{lem:seq}.

\section{Statement of the principal lemmas}\label{sec:lems}

We prove Lemma~\ref{lem:balanced} using the absorbing method of R\"odl, Ruci\'nski, and Szemer\'edi.  As is typical with this method, we have connecting, absorbing, and covering lemmas.

\begin{lemma}[Connecting lemma]\label{lem:connecting}
  For every $r \ge 2$ and $0 < \nu \le \frac{1}{r}$ there exists $\tau > 0$
  such that the following holds for every $n$.
  Let $G$ be an $r$-partite graph 
  with ordered partition $\cP = (V_1, \dotsc, V_r)$.
  Let $\ell = r(2r - 2)$.
  Suppose that $(U_1, \dotsc, U_r)$ is a sequence of sets 
  such that $U_i \subseteq V_i$ for $i \in [r]$, $U=\bigcup_{i=1}^r U_i$, and
  \begin{equation}\label{eq:con_min_deg}
    \text{for every $i \in [r]$ and $v \in V \setminus V_i$, 
      $|U_i| \ge \nu n$ and 
      $\deg_G(v, U_i) \ge \left(1 - \frac{1}{r} + \nu\right)|U_i|$}.
  \end{equation}
  Then for every pair of properly terminated $(r-1)$-walks $P_1$ and $P_2$ in $G$,
  there exist at least $\tau n^\ell$  $(r-1)$-walks $Q$ of length $\ell$ 
  contained in $U_1 \cup \dotsm \cup U_r$ 
  such that $P_1QP_2$ is a properly terminated $(r-1)$-walk.
\end{lemma}

\begin{lemma}[Absorbing lemma]\label{lem:absorbing}
  For $r \ge 2$, suppose that $\frac{1}{n} \ll \beta \ll \gamma < \frac{1}{r}$
  and let $G$ be a balanced $r$-partite graph on $n$ vertices
  with ordered partition $\cP = (V_1, \dotsc, V_r)$ such that 
  \begin{equation*}
    \delta_{\cP}(G) \ge 1 - \frac{1}{r} + \gamma.
  \end{equation*}
  Then there exists a properly terminated
  $(r-1)$-path $P_\textrm{abs}$ 
  such that $|V(P_{abs})| \le \beta n$,
  and, for every balanced set $Z \subseteq V(G) \setminus V(P_{abs})$ 
  for which $|Z| \le \beta^2 n$, 
  there exists a Hamiltonian $(r-1)$-path of $G[V(P_\textrm{abs}) \cup Z]$ that
  begins with the same $(r-1)$ vertices as $P_\textrm{abs}$
  and ends with the same $(r-1)$ vertices as $P_\textrm{abs}$.
\end{lemma}

\begin{lemma}[Covering lemma]\label{lem:covering}
  For $r \ge 2$, suppose that $\frac{1}{n} \ll \frac{1}{M_0} \ll \alpha \ll \gamma < \frac{1}{r}$
  and let $G$ be a balanced $r$-partite graph on $n$ vertices 
  with ordered partition $\cP = (V_1, \dotsc, V_r)$ and 
  \begin{equation*}
    \delta_{\cP}(G) \ge 1 - \frac{1}{r} + \gamma.
  \end{equation*}
  For some $M \le M_0$, there exist vertex
  disjoint properly terminated $(r-1)$-paths $P_1, \dotsc, P_M$ such that
  $W = V(G)\setminus \bigcup_{i=1}^{M} V(P_i)$ is balanced and $|W|\leq \alpha n$.
\end{lemma}

Before proving these three lemmas, we first show how to use Lemmas~\ref{lem:connecting}, \ref{lem:absorbing}, and \ref{lem:covering} to prove the balanced case of Theorem \ref{thm:main}.

\begin{proof}[Proof of Lemma~\ref{lem:balanced}]
  We can select $M_0$, $\nu$, $\alpha$, and $\beta$ so that 
  \begin{equation*}
    \frac{1}{n} \le \frac{1}{n_0} \ll \frac{1}{M_0} \ll \nu , \alpha \ll \beta \ll \gamma.
  \end{equation*}
  By Lemma~\ref{lem:absorbing}, 
  there exists a properly terminated-$(r-1)$ path 
  $P_\textrm{abs}$ disjoint from $K$ and $K'$ such that
  \begin{itemize}
    \item $|P_\textrm{abs}| \le \beta n$; and 
    \item for every balanced set $Z \subseteq V(G)$ such that $|Z| \le \beta^2 n$
      there exists a Hamiltonian $(r-1)$-path of $G[V(P_\textrm{abs}) \cup Z]$
      that starts and ends with the same $(r-1)$-vertices as $P_\textrm{abs}$.
  \end{itemize}
  Let $G' = G - \left(V(P_\textrm{abs}) \cup V(K) \cup V(K')\right)$.

  Uniformly at random select subsets $U_1, \dotsc, U_r$ such that
  for every $i \in [r]$, $U_i \subseteq V(G') \cap V_i$ and 
  $|U_i| = \ceiling{\nu n}$.
  By the Chernoff and union bounds, there exists an outcome such that
  \eqref{eq:con_min_deg} holds.
  Fix such an outcome and let 
  $U = U_1 \cup \dotsm \cup U_r$ and let $G'' = G' - V(U)$.

  By Lemma~\ref{lem:covering}, for some $M \le M_0$, there exist vertex
  disjoint properly terminated $(r-1)$-paths $P_1, \dotsc, P_M$ in $G''$ such that
  $W = V(G'')\setminus \bigcup_{i=1}^{M} V(P_i)$ is balanced and $|W| \leq \alpha n$.
  Since \eqref{eq:con_min_deg} holds, Fact~\ref{fact:walks_to_paths} 
  and Lemma~\ref{lem:connecting}
  imply that we can find $m+2$ disjoint $(r-1)$-paths, each of length
  $\ell = 2r(r-2)$, in $G[U]$ that connect 
  \begin{itemize}
    \item $K$ to $P_\textrm{abs}$, 
    \item $P_\textrm{abs}$ to $P_1$, 
    \item $P_{i-1}$ to $P_i$, for $2 \le i \le M$; and
    \item $P_M$ to $K'$
  \end{itemize}
  to form a $(r-1)$-path $P$.
  Let $Z = |V(G) \setminus V(P)|$,
  and note that 
  \begin{equation*}
    |Z| \le |U| + |W| \le r \ceiling{\nu n} + \alpha n \le \beta^2 n.
  \end{equation*}
  Therefore, there exists a Hamiltonian $(r-1)$-path of $G[P_\textrm{abs} \cup Z]$
  that starts and ends with the same $(r-1)$ vertices as $P_\textrm{abs}$.
  If $v_i = v'_i$ for every $i \in [r]$, then we have constructed a 
  Hamiltonian $(r-1)$-cycle.
  If $v_i \neq v'_i$ for every $i \in [r]$, then we have constructed a 
  Hamiltonian $(r-1)$-path that starts with $K$ and ends with $K'$.
\end{proof}

\section{Proof of the Connecting Lemma (Lemma~\ref{lem:connecting})}\label{sec:con}

Although we present the proof of Lemma~\ref{lem:connecting} in full, it closely follows proofs of similar lemmas given 
in \cite{jamshed2010embedding} and \cite{herdade2015stability}. 

\begin{definition}\label{def:poor}
  Let $G$ be a graph on $n$ vertices.
  For $U \subseteq V(G)$, we say that $W$ is \textit{$(U, \sigma)$-rich} if there are at least $\sigma n$
  vertices $u \in U$ for which $N(u)$ contains $W$, otherwise $W$ is called \textit{$(U,\sigma)$-poor}.  \end{definition}

The following simple observation and fact are critical for the inductive proof of the connecting lemma.
\begin{observation}\label{obs:rich}
  For $r \ge 3$, let $G$ be a graph on $n$ vertices,
  let $\cP = (V_1, \dotsc, V_r)$ be an ordered partition of $V(G)$, 
  and let $U_r \subseteq V_r$.
  Suppose that 
  \begin{equation*}
    W = x_1, \dotsc, x_{r-1}, z^1_1, \dotsc, z^1_{r-1}, \dotsc, z^s_1, \dotsc, z^s_{r-1}, y_1, \dotsc, y_{r-1}
  \end{equation*}
  is an $(r-2)$-walk of length $(s+2)(r-1)$ such that $W \cap V_r = \emptyset$ that is $(U_r, \sigma)$-rich.
  Then, by the definition of $(U_r, \sigma)$-rich, 
  there are at least $(\sigma n)^{s+1}$ tuples $(w^0, \dotsc, w^{s})$ such that
  $\{w^0, \dotsc, w^s\} \subseteq U_r$ and 
  $N(w^i)$ contains $W$ for each $0 \le i \le s$.  Therefore, for each such tuple
  \begin{equation*}
    x_1, \dotsc, x_{r-1}, w^0, z^1_1, \dotsc, z^1_{r-1}, w^1, \dotsc, z^s_1, \dotsc, z^s_{r-1}, w^{s}, y_1, \dotsc, y_{r-1}
  \end{equation*}
  is an $(r-1)$-walk of length $(s+2)r - 1$.
\end{observation}

By double counting, 
the following fact formalizes the observation that 
most neighborhoods do not contain many poor paths
for the simple reason that, by definition, 
poor paths are not contained in many neighborhoods,
\begin{fact}\label{fact:poor}
  For $r \ge 3$, $p \ge 0$ and $\sigma > 0$ the following holds.
  If $G$ is a graph on $n$ vertices and $U \subseteq V(G)$, 
  then there are at least $|U| -  \sigma n$ vertices $u \in U$ 
  such that only at most $\sigma n^p$ of the 
  $(r-2)$-walks of length $p$ contained in $N(u)$ are $(U,\sigma^2)$-poor.
\end{fact}
\begin{proof}
  Let $\cV_{poor}$ be the set of ordered $(p+1)$-tuples $(u, v_1, \dotsc, v_p) \in V^{p + 1}$
  such that 
  \begin{itemize}
    \item $u \in U$, 
    \item $W = v_1, \dotsc, v_p$ is a $(U, \sigma^2)$-poor $(r-2)$-walk, and
    \item $N(u)$ contains $W$. 
  \end{itemize}
  Because the number of ordered $p$-tuples is at most $n^p$,
  we have that $|\cV_{poor}| \le \sigma^2 n^{p + 1}$ (c.f.\ Definition~\ref{def:poor}).
  Let $U' \subseteq U$ be the set of vertices $u \in U$ such that 
  more than $\sigma n^p$ of the $(r-2)$-walks of length $p$ contained
  in $N(u)$ are $(U, \sigma^2)$-poor. Then, 
  \begin{equation*}
    |U'| \cdot \sigma n^p \le |\cV_{poor}| \le \sigma^2 n^{p + 1}.
  \end{equation*}
  Therefore, $|U'| \le \sigma n$ and the conclusion follows.
\end{proof}

\begin{proof}[Proof of Lemma~\ref{lem:connecting}]
  We will prove the lemma by induction on $r$.
  For the base case, note that when $r  = 2$,
  we have $\ell = 4$ and, by \eqref{eq:con_min_deg},
  the statement easily holds with $\tau = \nu^5 / 4$. 
  To see this, note that we can select 
  vertices $x_1, y_1 \in U_1$, and $y_2 \in U_2$ such that
  $P_1x_1$ and $y_1y_2P_2$ are $1$-paths. This can be done with $(1/2+\nu)|U_2|$ choices for $y_2$, $(1/2+\nu)|U_1|$ choices for $y_1$ and $(1/2+\nu)|U_1|$ choices for $x_1$ (recall that we only require $(r-1)$-walks). This gives at least $\left(\frac{\nu n}{2}\right)^3$ total selections and for every such selection
  we have
  \begin{equation*}
    |N(x_1) \cap N(y_1) \cap U_2|
    \ge \deg_G(x_1, U_2) + \deg_G(y_1, U_2) - |U_2| \ge 2\nu |U_2| \ge 2 \nu^2 n.
  \end{equation*}

  For the induction step, let $r \ge 3$ and suppose that the result holds for $r-1$.  Let $s = 2(r-1) - 2$, $q  = (r-1)(2(r-1) - 2) = (r-1)s$,
  and $p = q + 2(r - 1)$, and note that 
  \begin{equation}\label{eq:p+s+1}
    p + s + 2 = \left((r-1)s + 2(r - 1)\right) + s + 2 = r(s + 2) = 2r(r - 1) = \ell.
  \end{equation}
  Applying the induction hypothesis with $\nu/2$, $r-1$,  and $q$ playing the roles
  of $\nu$, $r$, and $\ell$ respectively we get that there exists $\mu > 0$ (playing the role of $\tau$) such that the following holds.
  \begin{claim}\label{clm:con_ind_hyp}
    If 
    $U'_i \subseteq U_i$  such that $|U'_i| \ge \nu n/2$ for all $i \in [r-1]$, and
    \begin{equation}\label{eq:con_min_deg_ind}
      \text{
        $\deg_G(v, U'_i) \ge \left(1 - \frac{1}{r-1} + \nu \right)|U'_i|$
      for all $v \in V \setminus V_i$}, 
    \end{equation}
    then for every pair of $(r-2)$-walks $x_1,\dotsc,x_{r-1}$ and 
    $y_1, \dotsc,y_{r-1}$
    such that $x_i,y_i \in U'_i$ for all $i \in [r-1]$ 
    there exist at least $\mu n^q$ $(r-2)$-walks of length $q$
    contained in $U'_1 \cup \dotsc \cup U'_{r-1}$ such that 
    $x_1,\dotsc,x_{r-1},Q',y_1,\dotsc,y_{r-1}$ is an $(r-2)$-walk.
  \end{claim}
 
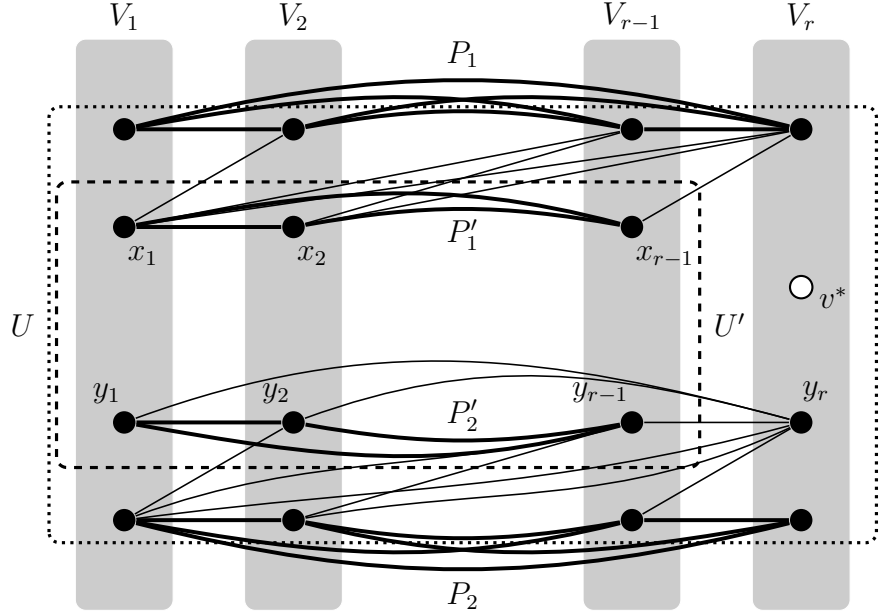
\begin{figure}
    \centering
    \def\xgap{2.25}
    \def\ygap{1.3}
    \def\wid{1.3}
    \begin{tikzpicture}
        \begin{pgfonlayer}{foreground}
            \foreach \i in {1,2,4,5}{
                \foreach \j in {-2,2}{
                    \node [vertex] (\j;\i) at (\i*\xgap,\j*\ygap) {};
                }
            }
            \foreach \i in {1,2}{
                \def\j{1};
                \node[vertex, label={[xshift=-6]-45:$x_{\i}$}] (\j;\i) at (\i*\xgap,\j*\ygap) {};
                \def\j{-1};
                \node[vertex, label={[xshift=6]135:$y_{\i}$}] (\j;\i) at (\i*\xgap,\j*\ygap) {};
            }
            \def\j{1};
            \node[vertex, label={[xshift=-6]-45:$x_{r-1}$}] (\j;4) at (4*\xgap,\j*\ygap) {};
            \def\j{-1};
            \node[vertex, label={[xshift=6]135:$y_{r-1}$}] (\j;4) at (4*\xgap,\j*\ygap) {};
            \node[vertex, label={[xshift=-6]60:$y_{r}$}] (-1;5) at (5*\xgap,-\ygap) {};
            \node[vertex, draw=black, fill=white,label={[xshift=-2,yshift=-4]0:$v^*$}] (v) at (5*\xgap,0.5) {};
        \end{pgfonlayer}
        \begin{pgfonlayer}{foreground}
            \node at (1*\xgap,2*\ygap+1.5) {$V_1$};
            \node at (2*\xgap,2*\ygap+1.5) {$V_2$};
            \node at (4*\xgap,2*\ygap+1.5) {$V_{r-1}$};
            \node at (5*\xgap,2*\ygap+1.5) {$V_r$};
            \node at (0.9,0) {$U$};
            \node at (4*\xgap+1.3,0) {$U'$};
            \node at (3*\xgap,2*\ygap+1.0) {$P_1$};
            \node at (3*\xgap,-2*\ygap-1.0) {$P_2$};
            \node at (3*\xgap,\ygap-0.1) {$P_1'$};
            \node at (3*\xgap,-\ygap+0.1) {$P_2'$};
        \end{pgfonlayer}
        \begin{pgfonlayer}{main}
            \foreach \j in {-2,2}{
                \draw [edge] (\j;1) -- (\j;2);
                \draw [edge] (\j;4) -- (\j;5);
                \draw [edge] (\j;2) to[out=\j*5, in=\j*85] (\j;4);
                \draw [edge] (\j;1) to[out=\j*5, in=\j*83] (\j;4);
                \draw [edge] (\j;2) to[out=\j*7, in=\j*85] (\j;5);
                \draw [edge] (\j;1) to[out=\j*7, in=\j*83] (\j;5);
            }
            \foreach \j in {-1,1}{
                \draw [edge] (\j;1) -- (\j;2);
                \draw [edge] (\j;2) to[out=\j*10, in=\j*170] (\j;4);
                \draw [edge] (\j;1) to[out=\j*10, in=\j*165] (\j;4);
            }
            \draw [edge, semithick] (1;1) -- (2;2);
            \draw [edge, semithick] (1;1) -- (2;4);
            \draw [edge, semithick] (1;1) -- (2;5);
            \draw [edge, semithick] (1;2) -- (2;4);
            \draw [edge, semithick] (1;2) -- (2;5);
            \draw [edge, semithick] (1;4) -- (2;5);
            \draw [edge, semithick] (-1;5) to[out=195,in=7.5] (-2;1);
            \draw [edge, semithick] (-1;5) to[out=205,in=10] (-2;2);
            \draw [edge, semithick] (-1;5) -- (-2;4);
            \draw [edge, semithick] (-1;4) to[out=195,in=20] (-2;1);
            \draw [edge, semithick] (-1;4) -- (-2;2);
            \draw [edge, semithick] (-1;2) -- (-2;1);
            \draw [edge, semithick] (-1;5) -- (-1;4);
            \draw [edge, semithick] (-1;5) to[out=165,in=20] (-1;2);
            \draw [edge, semithick] (-1;5) to[out=165,in=20] (-1;1);
        \end{pgfonlayer}
        \begin{pgfonlayer}{background}
            \draw [dotted,rounded corners, black, very thick] (1*\xgap-1.0,-2*\ygap-0.3) rectangle (5*\xgap+1.0,2*\ygap+0.3);
            \draw [dashed, rounded corners, black, very thick] (1*\xgap-0.9,-\ygap-0.6) rectangle (4*\xgap+0.9,\ygap+0.6);
        \end{pgfonlayer}
        \begin{pgfonlayer}{backbackground}
            \foreach \j in {1,2,4,5}{
                \draw [rounded corners, white, fill=gray!40] (\j*\xgap-\wid/2,-2*\ygap-1.2) rectangle (\j*\xgap+\wid/2,2*\ygap+1.2);
            }
        \end{pgfonlayer}
    \end{tikzpicture}
    \caption{Using induction to build the desired connection between $P_1$ and $P_2$ for Lemma~\ref{lem:connecting}.}
    \label{fig:connect}
\end{figure}

  Pick $\tau, \sigma > 0$ so that $\tau \ll \sigma \ll \mu, \nu$. 
  First note that, by \eqref{eq:con_min_deg}, there are at least $\frac{\nu n}{r} \ge \sigma n$ ways to select 
  $y_r \in U_r$ so that $y_rP_2$ is an $(r-1)$-path.
  Next, because $|U_r| \ge \nu n > \sigma n$, 
  Fact~\ref{fact:poor} implies that there exists $v^* \in U_r$ such that
  \begin{equation}\label{eq:vstar}
    \text{at most $\sigma n^{p}$ 
    of the $(r-2)$-walks of length $p$ contained in $N(v^*)$ are
    $(U_r, \sigma^2)$-poor.}
  \end{equation}
  For every $i \in [r-1]$, let $U'_i = N(v^*, U_i)$.
  Note that $|U'_i| \ge \frac{r-1}{r}|U_i| \ge \nu n/2$ and 
  for every $v \in V \setminus V_i$
  \begin{equation*}
    \deg(v, U'_i) \ge |U'_i| - \left(\frac{1}{r} - \nu\right)|U_i| \ge |U'_i| - \left(\frac{1}{r} - \nu\right)\frac{r}{r-1} |U'_i| \ge 
    \left(1 - \frac{1}{r-1} + \nu\right)|U'_i|. 
  \end{equation*} 
  Therefore, we can iteratively prepend vertices $y_{r-1}, \dotsc, y_1$ to $y_rP_2$ 
  and append vertices $x_1, \dotsc, x_{r-1}$ to $P_1$
  in at least $\left(\frac{\nu^2 n}{2}\right)^{2r-2}$ ways (Fig.~\ref{fig:connect}) so that the following holds:
  \begin{itemize}
    \item $x_i,y_i \in U'_i$ for $i \in [r-1]$; and
    \item both $P_1, x_1, \dotsc, x_{r-1}$ and
      $y_1, \dotsc, y_{r-1}, y_r, P_2$ are $(r-1)$-walks.
  \end{itemize}
  By Claim~\ref{clm:con_ind_hyp}, 
  the number of $(r-2)$-walks $Q'$ of length $q$ 
  contained in $U'_1 \cup \dotsm \cup U'_{r-1}$ such 
  that $x_1,\dotsc,x_{r-1},Q',y_1,\dotsc, y_{r-1}$ 
  is an $(r-2)$-path is at least $\mu {n}^q$.
   
  Therefore, there are at least
  \begin{equation*}
    \left(\frac{\nu^2 n}{2}\right)^{2r-2} \cdot \mu n^q = 2^{-2r + 2} \cdot \nu^2 \cdot \mu n^p \ge 2 \sigma n^{p}
  \end{equation*}
  $(r-1)$-walks
  \begin{equation*}
    x_1, \dotsc, x_{r-1}, Q', y_1, \dotsc, y_{r-1} =
    x_1, \dotsc, x_{r-1}, z^1_1, \dotsc, z^1_{r-1}, \dotsc, z^s_1, \dotsc, z^s_{r-1}, y_1, \dotsc, y_{r-1}
  \end{equation*}
  such that 
  \begin{itemize}
    \item $N(v^*)$ contains $x_1,\dotsc,x_{r-1},Q'y_1,\dotsc,y_{r-1}$;
    \item $x_1,\dotsc,x_{r-1},Q',y_1,\dotsc,y_{r-1}$ is an $(r-2)$-walk of length $p$; and
    \item both $P_1,x_1,\dotsc,x_{r-1}$ and $y_1,\dotsc,y_r,P_2$ are $(r-1)$-walks.
  \end{itemize}
  By \eqref{eq:vstar}, only $\sigma n^{p}$ of these paths are
  $(U_r, \sigma^2)$-poor so 
  at least $\sigma n^{p}$ of these paths are $(U_r, \sigma^2)$-rich.
  By Observation~\ref{obs:rich}, for every such $(U_r, \sigma^2)$-rich walk,
  there are at least $\left(\sigma^2 n\right)^{s+1}$ ordered tuples 
  $(w^0, \dotsc, w^{s})$ such that $\{w^0, \dotsc, w^s\} \subseteq U_r$ and
  \begin{equation*}
    x_1, x_2, \dotsc, x_{r-1}, w^0, z^1_1, \dotsc, z^1_{r-1}, w^1, \dotsc, z^s_1, \dotsc, z^s_{r-1}, w^{s}, y_1, \dotsc, y_{r-1}  
  \end{equation*}
  is an $(r-1)$-walk of length $p + s + 1= \ell-1$ (c.f.\ \eqref{eq:p+s+1}).
  Recalling that there were at least $\sigma n$ ways to select $y_r$
  gives us that 
  the number of $(r-1)$-walks $Q$ of length $\ell$ such that $P_1QP_2$ is an $(r-1)$-walk
  is at least 
  $\sigma n \cdot \sigma n^{p} \cdot \left(\sigma^2 n\right)^{s + 1} = \sigma^{2s + 4} n^{\ell} \ge \tau n^\ell$.
\end{proof}

\section{Proof of the Absorbing Lemma (Lemma~\ref{lem:absorbing})}\label{sec:abs}

\begin{definition}
  Let $2 \le r \le \ell$,
  let $G$ be an $r$-partite graph, and let $X$ be a balanced
  subset of $V(G)$. 
  A properly terminated $(r-1)$-path $a_1, \dotsc, a_\ell$ in $G$
  is an \textit{absorber of $X$}
  if there is an ordering of the vertices $\{a_1, \dotsc, a_\ell\} \cup X$ 
  that starts with the sequence $a_1, \dotsc, a_{r-1}$ 
  and ends with the sequence $a_{\ell - r + 1}, \dotsc, a_\ell$ 
  that is an $(r-1)$-path in $G$.
\end{definition}

The proof of the absorbing lemma follows by a standard probabilistic argument 
after the proof of the Lemma~\ref{lem:main_absorbing} below.

We will use the well known ``supersaturation'' result of Erd\H{o}s \cite{E} (see \cite[Theorem 2.11]{N}).

\begin{theorem}[Supersaturation]\label{thm:supersat}
 For all $r\geq 2$, $c'>0$, and positive integers $s_1, \dots, s_r$, there exists $n_0$ and $c$ such that if $G$ is a $r$-partite $r$-uniform hypergraph with ordered partition $(V_1, \dots, V_r)$ and at least $c'n^r$ edges, then $G$ contains at least $cn^{s_1+s_2+\dots+s_r}$ complete $r$-partite graphs with $s_i$ vertices in $V_i$ for all $i\in [r]$.
\end{theorem}

\begin{lemma}\label{lem:main_absorbing}
  For all $r \ge 2$ and $\frac{1}{n} \ll \alpha \ll \alpha' \ll \gamma \ll \frac{1}{r}$ the following holds
  with $\ell = 3r^2 - r$:
  
  Let $G$ be a balanced $r$-partite graph on 
  $n$ vertices with ordered partition $(V_1,\dotsc,V_r)$ 
  such that 
  $\delta_{\cP}(G) \ge 1 - \frac{1}{r} + \gamma$.
  If $X \subseteq V(G)$ is a balanced set of size $r$,
  then there are at least $(\alpha n)^{\ell}$ absorbers of $X$ in $G$.
\end{lemma}
\begin{proof}
  Let $x_1,\dotsc,x_r$ be an ordering of $X$ such that
  $x_i \in V_i$ for $i \in [r]$.
  
  We first describe what an absorber of $X$ will look like.  Suppose $$P=v_1^1\dots v_r^1 v_1^2\cdots v_r^2 \cdots v_1^{r-1} \cdots v_r^{r-1} v_1^r\cdots v_r^r$$ is an $(r-1)$-path of order $r^2$ where $v_i^j\in V_i$ for all $i\in [r]$.  For all $i,j\in [r]$, set $s_i^j=2$ if $i=j$ and $s_i^j=3$ otherwise.  
  
  Let $\cP$ be the $(s_1^1, \dots, s_r^1, \dots, s_1^r, \dots, s_r^r)$-blow up of $P$ where $D_i^j$ is the set corresponding to $v_i^j$. That is, replace each vertex $v_i^j$ with a set $D_i^j$ of order $s_i^j$, and if $\{v_i^j, v_{i'}^{j'}\}$ is an edge of $P$, add all edges between $D_{i}^{j}$ and $D_{i'}^{j'}$.
  
  We claim that, if we suppose that $D_1^i \cup \cdots \cup D_{i-1}^{i} \cup D_{i+1}^{i} \cup \cdots \cup D_r^{i}\subseteq N(x_i)$, for all $i \in [r]$, then $\cP$ contains an absorber of $X$.  For all $i\neq j\in [r]$, label the vertices of $D_i^j$ as $a_i^j, b_i^j, c_i^j$ and label the vertices of $D_i^i$ as $a_i^i$ and $c_i^i$.  
  Let
    $$ Q_1 = a_1^1 \cdots a_r^1 x_1 b_2^1 \cdots b_r^1 c_1^1 \cdots c_r^1 a_1^2 \cdots a_r^2 b_1^2 x_2 b_3^2 \cdots b_r^2 c_1^2 \cdots c_r^2 \cdots a_1^r \cdots a_r^r b_1^r \cdots b_{r-1}^{r} x_r c_1^r \cdots c_r^r $$ 
    and 
    $$ Q_2 = a_1^1 \cdots a_r^1 c_1^1 b_2^1 \cdots b_r^1 a_1^2 c_2^1 \cdots c_r^1 b_1^2 a_2^2 \cdots a_r^2 c_1^2 c_2^2 b_3^2 \cdots b_r^2 a_1^3 a_2^3 c_3^2 \cdots c_r^2 \cdots b_1^r \cdots b_{r-1}^{r} a_r^r c_1^r \cdots c_r^r, $$
    i.e., $Q_2 = T_1 \cdots T_r$ where
    \begin{align*}
    &T_1 = a_1^1 \cdots a_r^1 c_1^1 b_2^1 \cdots b_r^1 a_1^2 c_2^1 \cdots c_r^1,\\
    &T_i = b^i_1 \cdots b^i_{i-1}a^i_i\cdots a^i_r c^i_1\cdots c^i_ib^i_{i+1}\cdots b^i_r a^{i+1}_1\cdots a^{i+1}_i c^i_{i+1} \cdots c^i_r \, \text{ for $2 \le i \le r-1$, and }\\
    &T_r = b_1^r \cdots b_{r-1}^{r} a_r^r c_1^r \cdots c_r^r.
    \end{align*}
    Note that $Q_1$ and $Q_2$ are properly terminated $(r-1)$-paths which start with the same $r$ vertices and end with the same $r$ vertices, so $\cP$ contains an absorber for $X$.  See Figure \ref{fig:3absorb_cluster_mod}.
    
    \begin{example}
        In the case of $r=3$, the $2$-paths $Q_1$ and $Q_2$ are as follows: 
        \begin{align*}
            Q_1 &=   a_1^1 a_2^1 a_3^1 x_1 b_2^1 b_3^1 c_1^1 c_2^1 c_3^1 a_1^2 a_2^2 a_3^2 b_1^2 x_2 b_3^2 c_1^2 c_2^2 c_3^2 a_1^3 a_2^3 a_3^3 b_1^3 b_2^3 x_3 c_1^3 c_2^3 c_3^3 \\
            Q_2 &=   a_1^1 a_2^1 a_3^1 c_1^1 b_2^1 b_3^1 a_1^2 c_2^1 c_3^1 b_1^2 a_2^2 a_3^2 c_1^2 c_2^2 b_3^2 a_1^3 a_2^3 c_3^2 b_1^3 b_2^3 a_3^3 c_1^3 c_2^3 c_3^3 .
        \end{align*} \label{ex:r3}
    \end{example}
  
\begin{figure}
  \centering
    \begin{tikzpicture}
        \def\xgap{1.8}
        \def\xpt{0.2}
        \pgfmathsetmacro{\wid}{6*\xpt}
        \def\ypt{0.8}
        \pgfmathsetmacro{\yptD}{3*\ypt/2}
        \pgfmathsetmacro{\ypta}{\ypt/2}
        \pgfmathsetmacro{\yptb}{-\ypt/2}
        \pgfmathsetmacro{\yptc}{-3*\ypt/2}
        \pgfmathsetmacro{\yptE}{-2*\ypt}
        \pgfmathsetmacro{\yptX}{-3*\ypt}
        \pgfmathsetmacro{\hght}{4.2*\ypt}
        \newlength{\cent}
        \foreach \i in {1,2,3}{
            \foreach \j in {1,2,3}{ 
                \pgfmathtruncatemacro{\k}{3 * (\j - 1) + \i}
                \pgfmathsetmacro{\cent}{\k*\xgap}
                    \coordinate (C\k) at (\cent,0);
                    \coordinate (D\k) at ($(C\k)+(0,\yptD)$); 
                    \coordinate (E\k) at ($(C\k)+(0,\yptE)$); 
                \begin{pgfonlayer}{background}
                    \draw [rounded corners, fill=white, thick] (\k*\xgap-0.5*\wid,-0.5*\hght) rectangle (\k*\xgap+0.5*\wid,0.5*\hght);
                \end{pgfonlayer}
                \begin{pgfonlayer}{foreground}
                    \node at (D\k) {$\mathbf{D_\i^\j}$};
                \end{pgfonlayer}
            }
        }
        \begin{pgfonlayer}{backbackground}
            \foreach \k in {1, ..., 7} {
                \pgfmathtruncatemacro{\l}{\k + 1}
                \pgfmathtruncatemacro{\m}{\k + 2}
                \draw [bigedge, bend left=100] (D\k) to (D\l);
                \draw [bigedge, bend left=100] (D\k) to (D\m);
            }
            \draw [bigedge, bend left=100] (D8) to (D9);
        \end{pgfonlayer}
        \begin{pgfonlayer}{foreground}
            \foreach \j in {1,2,3}{
                \pgfmathtruncatemacro{\i}{1}
                \pgfmathtruncatemacro{\k}{3 * (\j - 1) + \i}
                \node [vertex, fill=black, label=left:{$a_\i^\j$}] (a\i\j) at ( $ (C\k) + (\xpt,\ypta) $ )  { };
                \ifthenelse{\i=\j}{}{\node [vertex, fill=black, label=left:{$b_\i^\j$}] (b\i\j) at ( $ (C\k) + (\xpt,\yptb) $ )  { }};
                \node [vertex, fill=black, label=left:{$c_\i^\j$}] (c\i\j) at ( $ (C\k) + (\xpt,\yptc) $ )  { };
                \pgfmathtruncatemacro{\i}{2}
                \pgfmathtruncatemacro{\k}{3 * (\j - 1) + \i}
                \node [vertex, draw=black, fill=gray!30, label=left:{$a_\i^\j$}] (a\i\j) at ( $ (C\k) + (\xpt,\ypta) $ )  { };
                \ifthenelse{\i=\j}{}{\node [vertex, draw=black, fill=gray!30, label=left:{$b_\i^\j$}] (b\i\j) at ( $ (C\k) + (\xpt,\yptb) $ )  { }};
                \node [vertex, draw=black, fill=gray!30, label=left:{$c_\i^\j$}] (c\i\j) at ( $ (C\k) + (\xpt,\yptc) $ )  { };
                \pgfmathtruncatemacro{\i}{3}
                \pgfmathtruncatemacro{\k}{3 * (\j - 1) + \i}
                \node [vertex, fill=gray!70, label=left:{$a_\i^\j$}] (a\i\j) at ( $ (C\k) + (\xpt,\ypta) $ )  { };
                \ifthenelse{\i=\j}{}{\node [vertex, fill=gray!70, label=left:{$b_\i^\j$}] (b\i\j) at ( $ (C\k) + (\xpt,\yptb) $ )  { }};
                \node [vertex, fill=gray!70, label=left:{$c_\i^\j$}] (c\i\j) at ( $ (C\k) + (\xpt,\yptc) $ )  { };
            }
            \node [vertex,fill=black,label=below:{$x_1$}] (x1) at ($(C1)+(\xpt,\yptX)$) {};
            \node [vertex,draw=black,fill=gray!30,label=below:{$x_2$}] (x2) at ($(C5)+(\xpt,\yptX)$) {};
            \node [vertex,fill=gray!70,label=below:{$x_3$}] (x3) at ($(C9)+(\xpt,\yptX)$) {};
        \end{pgfonlayer}
        \begin{pgfonlayer}{main}
            \draw [bend left=5] (c11) to (a21);
            \draw [bend left=24] (c11) to (a31);
            \draw [bend left=0] (a12) to (b21);
            \draw [bend left=15] (a12) to (b31);
            \draw [bend left=0] (b12) to (c21);
            \draw [bend left=15] (b12) to (c31);
            \draw [bend right=0] (c12) to (a22);
            \draw [bend right=12] (c12) to (a32);
            \draw [bend left=1] (c22) to (a32);
            \draw [bend right=12] (b32) to (a13);
            \draw [bend right=2] (b32) to (a23);
            
            \draw [bend right=14] (c32) to (b13);
            \draw [bend right=2] (c32) to (b23);
            \draw [bend left=15] (a33) to (c13);
            \draw [bend left=0] (a33) to (c23);
        \end{pgfonlayer}
        \begin{pgfonlayer}{backbackground}
            \draw [bigedge, bend right=10] (x1) to (E2);
            \draw [bigedge, bend right=10] (x1) to (E3);
            \draw [bigedge, bend left=10] (x2) to (E4);
            \draw [bigedge, bend right=10] (x2) to (E6);
            \draw [bigedge, bend left=10] (x3) to (E7);
            \draw [bigedge, bend left=10] (x3) to (E8);
        \end{pgfonlayer}
    \end{tikzpicture}
  \caption{An absorber for $X = \{x_1, x_2, x_3\}$ from Lemma~\ref{lem:main_absorbing}.  The edges between $D_i^j$'s and between $X$ and $D_i^j$ are indicated by solid black lines.  The edges of $Q_1$ are not shown.  The edges of $Q_2$ that are not in $Q_1$ are shown.}
  \label{fig:3absorb_cluster_mod}
\end{figure}
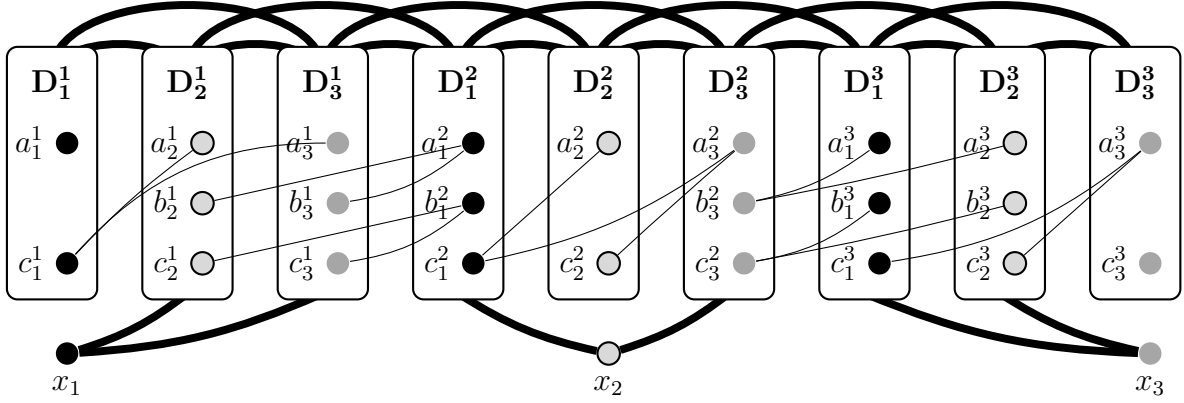
  
  Now we show that there are $\Omega(n^{3r^2-r})$ copies of $\cP$ which contain the absorber of $X$ as described above.  By a Chernoff bound (Theorem~\ref{chernoff}), for all $i\in [r]$ there exists a partition $V_i=\{V_i^1, \dots, V_i^r\}$ such that for all $i,j\in [r]$ and all $v\in V(G)\setminus V_i$, $$\deg\left(v, V_i^j\right)\geq \left(1-\frac{1}{r}+\frac{\gamma}{2}\right)\left|V_i^j\right|.$$ 
  
  One can see that constructing greedily (from the middle out), there are at least $(\frac{\gamma}{2} n)^{r^2}$ properly ordered $(r-1)$-paths $P=v_1\dots v_{r^2}$ of order $r^2$ such that for all $i\in [r]$, $$\{v_{ir+1}, \dots, v_{ir+i-1}, v_{ir+i+1}, \dots, v_{(i+1)r}\}\subseteq N(x_i).$$  Treating each such copy as an edge in an $r^2$-partite $r^2$-uniform hypergraph $H$ with ordered partition $(V_1^1, \dots, V_r^1, \dots, V_1^r, \ldots, V_r^r)$ and applying Theorem \ref{thm:supersat} to $H$, we have that there exists at least $\alpha n^{3r^2-r}$ copies of the $(s_1^1, \dots, s_r^1, \dots, s_1^r, \dots, s_r^r)$-blow up of $P$. 
\end{proof}

\begin{proof}[Proof of Lemma~\ref{lem:absorbing}]
  Let $\alpha$ be such that $\frac{1}{n} \ll \alpha \ll \beta$, let
  $\ell = 3r^2 - r$, and 
  let $\cA'$ be the collection of all ordered sequences 
  $(a_1, \dotsc, a_\ell)$ of vertices
  such that for every $i \in [\ell]$ and $j \in [r]$, if $a_i \in V_j$, then $i \equiv j \pmod r$.
  Let $\cX$ be the collection of all balanced $r$-subsets of $V(G)$.
  For every $X \in \cX$, let 
  \begin{equation*}
    \cA'_X = \{(a_1, \dotsc, a_\ell) \in \cA' : \text{$a_1, \dotsc, a_\ell$ is an absorber of $X$} \},
  \end{equation*}
  and note that, by Lemma~\ref{lem:main_absorbing}, we have 
  \begin{equation}\label{eq:size_of_AU}
    |\cA'_X| \ge (\alpha n)^\ell.
  \end{equation}

  Now create a random set $\cA_\textrm{ran}$ 
  by select each sequence in $\cA'$ independently at random with probability 
  $\rho = \beta^{1.1} n^{-\ell + 1}$, so since $|\cA'| = n^\ell$,
  \begin{equation*}
    \E |\cA_\textrm{ran}| = \rho |\cA'| \le \frac{\beta n}{4 \ell}, 
  \end{equation*}
  and, by \eqref{eq:size_of_AU}, for every $X \in \cX$,
  \begin{equation*}
    \E |\cA_{X}' \cap \cA_\textrm{ran}| \ge \rho (\alpha n)^{\ell} \ge 4 \beta^2 n.
  \end{equation*}
  So, by the Chernoff bound and the union bound, with high probability
  \begin{equation*}
    |\cA_\textrm{ran}| \le \frac{\beta n}{3 \ell} \qquad \text{ and } \qquad 
    |\cA_{X}' \cap \cA_\textrm{ran}| \ge 3 \beta^2 n \qquad \text{ for every $X \in \cX$}.
  \end{equation*}
  Let $\cA_\textrm{rep}$ contain the pairs of tuples in $\cA'$ in which
  a vertex is repeated, i.e.,
  \begin{equation*}
    \cA_\textrm{rep} = \{ \{S,T\} : S,T \in \cA', S \neq T, \text{ and a vertex appears at least twice in sequence $S,T$} \}.
  \end{equation*}
  We can construct every pair in $\cA_\textrm{rep}$ by selecting an arbitrary vertex,
  placing that vertex in $2$ of the $2 \ell$ possible entries, and then arbitrarily filling the
  remaining $2\ell - 2$ entries, so
  \begin{equation*}
    \E |\cA_\textrm{rep} \cap \cA_\textrm{ran}| =
    \rho^2 |\cA_\textrm{rep}| \le
    \rho^2 \cdot n \cdot \binom{2\ell}{2} \cdot n^{2\ell - 2} \le \beta^2 n.
  \end{equation*}
  By the Markov bound, with probability $1/2$, we have that $|\cA_\textrm{rep}| \le 2 \beta^2 n$.
  Therefore, there must exist some random outcome $\cA_\text{ran}$ such that
  if we remove every pair in $\cA_\textrm{rep} \cap \cA_\text{ran}$ and every sequence that
  is not absorbing for some $X \in \cX$ to form $\cA$ then we have that
  \begin{itemize}
    \item $|\cA| \le \beta n/(3 \ell)$;
    \item $|\cA \cap \cA_X'| \ge \beta^2 n$ for every $X \in \cX$;
    \item the sequences in $\cA$ are pairwise vertex-disjoint; and
    \item for every $P \in \cA$, $P$ is an absorber for some $X \in \cX$,
      so $P$ is an $(r-1)$-path.
  \end{itemize}
  
  Lemma~\ref{lem:connecting} (with $(V_1, \dotsc, V_r)$ and $\gamma$
  playing the roles of $(U_1, \dotsc, U_r)$ and $\nu$, respectively)
  and Fact~\ref{fact:walks_to_paths} together imply that we can
  connect the $(r-1)$-paths in $\cA$ (in an arbitrary order) 
  with paths of length $r(2r - 2) < 2 \ell$ 
  to form the desired absorbing $(r-1)$-path $P_\textrm{abs}$.
  We have that $|V(P_\textrm{abs})| = \ell |\mathcal{A}| + r(2r-2) (|\mathcal{A}| - 1) < 3\ell|\mathcal{A}| \le \beta n$.

  Let $Z \subseteq V(G) \setminus V(P_\textrm{abs})$ be a balanced set where $|Z| \le \beta^2 n$.
  We can partition $Z$ into balanced $r$-subsets so that each part is in $\cX$.  
  Since there are at most $|Z|/r < \beta^2 n$ parts in such a partition,
  we can greedily match each part $X$ to some path $P \in \cA \cap \cA_X$.
  Since $P$ is an absorber of $X$, we can construct the desired Hamiltonian
  $(r-1)$-path of $G[V(P_\textrm{abs}) \cup Z]$.
\end{proof}

\section{The regularity lemma}\label{sec:reg}

We now review Szemer\'edi's well-known regularity lemma \cite{Sz}.

\begin{definition}
In a graph $G$, for each pair of disjoint non-empty sets $A, B \subseteq V(G)$ we write $G[A, B]$ for the bipartite subgraph of $G$ with vertex classes $A$ and $B$ and whose edges are all edges of $G$ with one endvertex in $A$ and the other in $B$, and denote the \emph{density} of $G[A, B]$ by $d_G(A, B) = \tfrac{e(G[A, B])}{|A||B|}$. 

We say that $G[A, B]$ is $(d,\eps)$-\emph{regular} if $d_G(X, Y) = d \pm \eps$ for every $X \subseteq A$ and $Y \subseteq B$ with $|X| \geq \eps |A|$ and $|Y| \geq \eps |B|$, and we write that $G[A, B]$ is $(\geq\!\!d, \eps)$-\emph{regular} to mean that $G[A, B]$ is $(d', \eps)$-regular for some $d' \geq d$. 

Also, we say that $G[A, B]$ is $(d, \eps)$-\emph{super-regular} if $G[A, B]$ is $(\geq\!\!d, \eps)$-regular, every vertex of $A$ has at least $(d-\eps)|B|$ neighbours in $B$, and every vertex of $B$ has at least $(d-\eps)|A|$ neighbours in $A$. 
\end{definition}

The following results are well-known elementary consequences of the definitions.
\begin{lemma}[Slicing Lemma]\label{lem:slicing}
For every $d, \eps, \beta > 0$, if $G[A, B]$ is $(d, \eps)$-regular, and $X \subseteq A$ and $Y \subseteq B$ have sizes $|X| \geq \beta |A|$ and $|Y| \geq \beta |B|$, then $G[X, Y]$ is $(d, \eps/\beta)$-regular.
\qed
\end{lemma} 

\begin{lemma}\label{lem:makesuperreg}
For every $d, \eps > 0$ with $\eps < \frac{1}{2}$, if $G[A, B]$ is $(\geq\!\!d, \eps)$-regular, then there are sets $X \subseteq A$ and $Y \subseteq B$ with sizes $|X| \geq (1-\eps) |A|$,  and $|Y| \geq (1-\eps) |B|$ such that $G[X, Y]$ is $(d, 2\eps)$-super-regular.
\qed
\end{lemma}

\begin{definition}
  Let $G$ be a graph on $n$ vertices and suppose that $\cC$ is a collection of disjoint subsets of $V(G)$.
  Define the \textit{$(G, \cC, d, \eps)$-cluster graph} to be the graph with vertex set $\cC$
  in which distinct $A,B \in \cC$ form an edge if $G[A,B]$ is $(\geq\!\!d, \eps)$-regular.
\end{definition}

\begin{definition}
  Let $\cP = (V_1, \dotsc, V_r)$ be an ordered partition of $V(G)$.
  We say that a collection $\cC$ of vertex disjoint subsets of $V(G)$
  \textit{respects $\cP$} if for every $C \in \cC$
  we have $C \subseteq V_i$ for some $i \in [r]$.
  If $\cC$ respects $\cP$, we let 
  $\cP(\cC)$ be the partition $(\cC_1, \dotsc, \cC_r)$ of
  $\cC$ in which every $C \in \cC$ is in
  $\cC_i$ when $C \subseteq V_i$.
\end{definition}

We now state the standard degree form of the regularity lemma.
\begin{lemma}[Degree Form of Szemer\'edi's Regularity Lemma]\label{lem:dreg}
  For every $\varepsilon > 0$ and $0 < d <1$ and integers
  $r$ and $N_0$ there exists $N_1$ such that the following holds.
  If $G$ is an $r$-partite graph on $n$ vertices 
  with ordered partition $\cP$, then
  there exists a partition $U_0,\dotsc, U_N$ of $V(G)$ 
  and a spanning subgraph $R$ of $G$ such that the following holds:
  \begin{itemize}
    \item $N_0 \le N \le N_1$;
    \item $|U_0| \le \eps n$;
    \item $|U_1| = \dotsm = |U_N|$;
    \item the collection $U_1, \dotsc, U_N$ respects the partition 
      $(V_1, \dotsc, V_r)$;
    \item $\deg_R(v) \ge \deg_G(v) - (d+\varepsilon)n$ for every $v \in V(G)$;
    \item $|E(R[U_i]) = 0$ for every $1 \le i \le N$; and
    \item for every $1 \le i < j \le N$, the graph $R[U_i,U_j]$ either
      $(\geq\!\!d, \eps)$-regular or has no edges.
  \end{itemize}
\end{lemma}

From the degree form of the regularity lemma, it is easy
to show that we have Lemma~\ref{lem:reg} below.
Since the proof is standard, we only provide a sketch.
\begin{lemma}\label{lem:reg}
  Suppose that
  \begin{equation*}
    \frac{1}{n} \ll \frac{1}{N_1} \ll \eps \ll d \ll \eta, \frac{1}{N_0}, \frac{1}{r}.
  \end{equation*}
  Let $G$ be a balanced $r$-partite graph on $n$ vertices 
  with ordered partition $\cP$.
  Then there exists $\cC$, which is a collection of vertex disjoint subsets of $V(G)$ and $R$ a spanning
  subgraph of $G$ such that 
  \begin{enumerate}[label=(R\arabic*)]
    \item\label{R1} $N_0 \le |\cC| \le N_1$;
    \item\label{R2} $\cC$ covers all but at most $\eps n$ vertices of $G$;
    \item\label{R3} every element in $\cC$ has the same order;
    \item\label{R4} $\cC$ respects the partition $\cP$ and the partition $\cP(\cC) = (\cC_1, \dotsc, \cC_r)$ is balanced;
    \item\label{R5} for every $v \in V(G)$, we have $\deg_R(v) \ge \deg_G(v) - (d + \eps)n$; 
    \item\label{R6} 
      for every $U \in \cC$, we have $E_R(U) = \emptyset$, and
      for every pair of distinct $A,B \in \cC$, either 
      $E(R[A, B]) = \emptyset$ or $R[A, B]$ is $(\ge\!\!d, \eps)$-regular; and
    \item\label{R7} if $\cG$ is the $(G, \cC, d, \eps)$-cluster graph, 
      then $\delta_{\cP(\cC)}(\cG) \ge \delta_{\cP}(G) - \eta$.
  \end{enumerate}
\end{lemma}
\begin{proof}[Proof sketch]
  Pick $\varepsilon'$ and $d'$ such that 
  $\frac{1}{N_1} \ll \varepsilon' \ll \varepsilon \ll d' \ll d$.
  Lemma~\ref{lem:dreg} 
  implies that
  there exists 
  a spanning subgraph $R$ of $G$
  and $U_0, U_1, \dotsc, U_N$ 
  a collection of vertex disjoint subsets of $V(G)$ such that
  the conclusions of Lemma~\ref{lem:dreg} hold
  with 
  $\varepsilon'$, $d'$, $r$ and $2N_0$ playing the roles of 
  $\varepsilon$, $d$, $r$ and $N_0$.
  In particular, we have that $N \ge 2 N_0$ and 
  $U_1, \dotsc, U_N$ covers all but at most 
  $\varepsilon' n$ of the vertices of $G$.
  Therefore, by removing a small fraction of the sets from
  the collection $U_1, \dotsc, U_N$ we can create $\cC$ 
  a collection of vertex disjoint subsets of $V(G)$ 
  such that \ref{R1} - \ref{R6} all hold.
  
  To see that \ref{R7} holds as well,
  let $\cP(\cC) = (\cC_1, \dotsc, \cC_r)$ and 
  let $i,j \in [r]$ such that $i \neq j$.
  For every $C \in \cC_i$ and $v \in C$, 
  \ref{R2}, \ref{R3}, \ref{R5}, and \ref{R6} imply that
  \begin{align*}
    \frac{\deg_{\cG}(C, \cC_j)}{|\cC_j|} \ge
    \frac{\deg_R(v, V(\cC_j))}{|C||\cC_j|} &\ge 
    \frac{\deg_R(v, V_j) - \varepsilon n}{n/r} \\
    &\ge \frac{\deg_G(v, V_j) - (d + \varepsilon)n - \varepsilon n}{n/r} \ge
    \delta_{\cP}(G) - \eta. \qedhere
  \end{align*}
\end{proof}

We make the following definition to help describe the version
of the well-known blow-up lemma that we will need.
\begin{definition}
  For a graph $R$ and $\cC$ be a collection of vertex disjoint 
  subsets of $V(R)$,
  we let $K(\cC, R)$ be the graph on $V(\cC)$ such
  that for every distinct $x,y \in V(R)$ the graph 
  $K(\cC, R)$ has the edge $\{x,y\}$ if and only if
  $x$ and $y$ are in distinct sets $A,B \in \cC$ and 
  $E(R[A,B]) \neq \emptyset$.

  For a subgraph $H$ of $K(\cC, R)$, 
  a \textit{copy of $H$ in $R$ that respects $\cC$}
  is an injective function $f : V(H) \to V(R)$ such
  that $\{x,y\} \in E(H)$ implies $\{f(x), f(y)\} \in E(R)$
  and, for every $v \in V(H)$ and $C \in \cC$,
  $v \in C$ implies $f(v) \in C$.
\end{definition}

\begin{lemma}[Blow-up Lemma \cite{KSS1}]\label{lem:blowup}
  Suppose that $\frac{1}{m} \ll \eps \ll d, \frac{1}{D}$.
  Let $G$ be a graph on $n$ vertices;
  let $\cC$ be a collection of vertex disjoint subsets of $V(G)$ each 
  of size $m$; and 
  let $R$ be a spanning subgraph of $G$ such that 
  for every $U \in \cC$, we have $E_R(U) = \emptyset$, and
  for every pair of distinct $A,B \in \cC$, either 
  $E(R[A, B]) = \emptyset$ or $R[A, B]$ is $(\ge\!\!dd, \eps)$-super-regular.
  If $H \subseteq K(\cC, R)$ and $\Delta(H) \le D$, 
  then there exists a copy of $H$ in $R$
  that respects $\cC$.
\end{lemma}

\section{Proof of the Covering Lemma (Lemma~\ref{lem:covering})}\label{sec:cov}

\begin{definition}
  Let $G$ be a graph and let $\cK$ be the copies of $K_r$ in $G$.
  A \textit{fractional $K_r$-tiling of a graph $G$} is a weight function $w : E(\cK) \to \mathbb{R}_{\ge 0}$ in which, for every $v \in V(G)$, 
  the sum of the weights on the copies of $K_r$ that contain $v$ is at most one.
  That is, we have that
  \begin{equation*}
    \sum \{w(K) : \text{$K \in \cK$ and $K$ contains $v$}\}  \le 1 \qquad \text{for every $v \in V(G)$}.
  \end{equation*}
  The \textit{size} of $w$ is $\sum \{w(K): K \in \cK \}$, 
  and we say that $w$ is \textit{perfect} if the size of $w$ is exactly $|V(G)|/r$.
  Note that $w$ is perfect if and only if 
  \begin{equation*}
    \sum \{w(K) : \text{$K \in \cK$ and $K$ contains $v$}\} = 1  \qquad \text{for every $v \in V(G)$}.
  \end{equation*}
\end{definition}

We will use the following lemma which can be found as a corollary to~\cite[Lemma 2.2]{MMS}. (See also,~\cite{LM1,KM1}.)
\begin{lemma}\label{lem:balanced_fractional}
  If $G$ is a balanced $r$-partite graph on $n$ vertices with partition 
  $\cP$ and $\delta_{\cP}(G) \ge 1 - \frac{1}{r}$, 
  then $G$ has a perfect fractional $K_r$-tiling.
\end{lemma}

\begin{remark}
  Here we could have shortened our proof by using existing results on 
  perfect $K_r$-tilings in multipartite graphs (see \cite{KM1, LM2}).
  We chose to only use the above lemma on 
  perfect fractional $K_r$-tilings, which is relatively short, 
  to make this paper more self-contained.
\end{remark}

The following lemma is a consequence of 
Lemma~\ref{lem:slicing} (The Slicing Lemma), 
Lemma~\ref{lem:makesuperreg}, and
Lemma~\ref{lem:blowup} (The Blow-up Lemma).
\begin{lemma}\label{lem:r-1paths}
  Let $\frac{1}{m} \ll \varepsilon \ll d \ll \alpha' < \frac{1}{r}$,
  let $G$ be an $r$-partite graph with ordered partition $(V_1, \dotsc, V_r)$ 
  and for $i \in [r]$, let $C_i$ be an $m$-subset of $V_i$.
  Suppose that the sets $C_1, \dotsc, C_r$ are pairwise $(\geq\!\!d, \varepsilon)$-regular
  and for every $i \in [r]$, we have $C'_i \subseteq C_i$.
  If $z$ is a positive integer such that $|C_i \setminus C_i'| + z \le (1 - \alpha')m$ for every $i \in [r]$,
  then there exists $P$ a properly terminated $(r-1)$-path in
  $G[C'_1 \cup \dotsm \cup C'_r]$ such that for every $i \in [r]$ 
  the path $P$ intersects $C'_i$ in exactly $z$ vertices.
\end{lemma}
\begin{proof}
  Note that the conditions imply that $|C'_i| \ge \alpha' m + z$ for every $i \in [r]$.
  So, Lemma~\ref{lem:slicing} (the Slicing Lemma), 
  implies that the sets $C'_1, \dotsc, C'_r$ are pairwise $(\geq\!\!d, \varepsilon^{2/3})$-regular.
  By applying Lemma~\ref{lem:makesuperreg} $\binom{r}{2}$ times,
  we can construct $C''_i \subseteq C'_i$ for $i \in [r]$ 
  such that $|C''_i| \ge z$
  and the sets $C''_1, \dotsc, C''_r$ 
  are pairwise $(d, \varepsilon^{1/3})$-super-regular.
  Lemma~\ref{lem:blowup} (the Blow-up Lemma)
  then implies the existence of the desired $(r-1)$-path $P$.
\end{proof}

\begin{proof}[Proof of Lemma~\ref{lem:covering}]
  Select constants $N_0$, $M_0$, $\eps$, $\alpha'$, $\eta$ and $d$ so that 
  \begin{equation*}
    \frac{1}{n} \ll \frac{1}{M_0} \ll \frac{1}{N_1} \ll \eps \ll \alpha' \ll \alpha \ll  d \ll \eta \ll \gamma < \frac{1}{r}. 
  \end{equation*}
  
  Lemma~\ref{lem:reg} implies the existence of a collection $\cC$ of disjoint subsets of $V(G)$ such that 
  \begin{itemize}
    \item $|\cC| \le N_1$
    \item $\cC$ covers all but at most $\eps n$ of the
      vertices in $V(G)$;
    \item there exists $m$ such that for every
      $C \in \cC$ we have $|C| = m$;
    \item $\cC$ respects $\cP$
      and if we let $\cP' = \cP(\cC)$
      and $\cG = (G, \cC, d, \eps)$,
      then $\cP'$ is balanced and 
      \begin{equation*}
        \delta_{\cP'}(G) \ge 1 - \frac{1}{r} + \frac{\gamma}{2}.
      \end{equation*}
  \end{itemize}

  Lemma~\ref{lem:balanced_fractional} implies that there 
  exists a perfect fractional $K_r$-tiling of $\cG$, 
  and let $\cK_1, \dotsc, \cK_M$ be an arbitrary
  ordering of the copies of $K_r$ in $\cG$ that receive positive
  weight in such a fractional $K_r$-tiling.
  Note that $M \le \binom{N_1}{r}$
  and that there are positive weights $w_1, \dotsc, w_M$ 
  such that for every $C \in \mathcal{C}$,
  \begin{equation*}
    \sum_{i=1}^{M} \left\{w_i : \text{$\cK_i$ contains the cluster $C$}\right\} = 1,
  \end{equation*}
  and $\sum_{i=1}^{M} w_i = |\cC|/r \ge (1 - \varepsilon)n/(mr)$.
  For each $i \in [M]$, let $z_i = \floor{(1-\alpha') w_i m}$
  and note that
  \begin{equation*}
    \sum_{i=1}^{M} z_i = \sum_{i=1}^{M} \floor{(1-\alpha')w_im} \ge (1-\alpha')|\cC|\frac{m}{r} - M \ge
    (1-\alpha')(1 - \varepsilon)\frac{n}{r} - M \ge \left(1 - \alpha\right)\frac{n}{r}.
  \end{equation*}

  We can now prove the lemma by constructing disjoint properly terminated $(r-1)$-paths $P_1, \dotsc, P_M$ such that 
  for each $i \in [M]$, the $(r-1)$-path $P_i$ has length exactly $r z_i$
  because then
    $\left|\bigcup_{i=1}^{m} V(P_i)\right| \ge (1 - \alpha)n$.

  To see that such a construction is possible, assume that,
  for some $t \in [M]$, we have constructed $t-1$ disjoint properly terminated $(r-1)$-paths
  $P_1, \dotsc, P_{t-1}$ such that for every $j \in [t-1]$
  the path $P_j$ is contained in the clusters of $\cK_j$
  and for every cluster $C$ contained in $\cK_j$ the $(r-1)$-path $P_j$ 
  intersects $C$ in exactly $z_j$ vertices.
  
  Let $C_1, \dotsc, C_r$ be the clusters in $\cK_t$.
  We can assume that $C_i \subseteq V_i$ for $i \in [r]$ since the partition $\cC$ respects the partition $\cP$
  and the clusters $C_1, \dotsc, C_r$ are pairwise $(\geq\!\! d, \varepsilon)$-regular.
  For $i \in [r]$, let $C'_i \subseteq C_i$ be the vertices in $C_i$ that do not 
  intersect one of the previously constructed paths $P_1, \dotsc, P_{t-1}$.
  Recall that for $i \in [r]$, we have that $\cK_t$ contains the cluster $C_i$, so 
  \begin{align*}
    |C_i \setminus C'_i| + z_t &= 
    \left(\sum_{j=1}^{t-1} \left\{ z_j : \text{$\cK_j$ contains the cluster $C_i$} \right\}\right) + z_t \\
    &\le \sum_{j=1}^{M} \left\{ z_j : \text{$\cK_j$ contains the cluster $C_i$} \right\} \\
    &\le \sum_{j=1}^{M} \left\{ (1-\alpha')w_j m : \text{$\cK_j$ contains the cluster $C_i$} \right\} = (1 - \alpha') m.
  \end{align*}
  Therefore, Lemma~\ref{lem:r-1paths} implies that there exists 
  an $(r-1)$-path $P_t$ contained in $G[C'_1 \cup \dotsm \cup C'_r]$ 
  such that, for $i \in [r]$, the path $P_t$ intersects $C'_i$ in exactly $z_t$ vertices.
\end{proof}

\section{Proof of the Partitioning and Sequencing Lemma (Lemma \ref{lem:seq})}\label{sec:seq}
Before we begin the proof, we give some further terminology and observations regarding properly ordered paths.

For $1 \le k' \le k$, we say that an $(r-1)$-path $v_1, \dotsc, v_{k'}$ 
is \textit{increasing} if for every $1 \le i < i' \le k'$ we have that $v_i \in V_j$ and
$v_{i'} \in V_{j'}$ with $j < j'$, i.e., a path is increasing if it traverses the
sets $V_1, \dotsc, V_k$ in order (though it might skip any number of the sets).
All of the paths that we will construct can be partitioned into subpaths
on either $r$ or $r+1$ vertices that are increasing.
We call such an $(r-1)$-path \textit{properly ordered}. 
We now give a more formal definition.
\begin{definition}[Properly ordered/$j$-th subsequence]
Let $P=v_1v_2\cdots v_p$ be an $(r-1)$-path and let $f:[p]\rightarrow [k]$ be such that $v_{f(i)}\in V_j$.  We say that $P$ is \textit{properly ordered} if there exists $0=p_0, p_1, \dots, p_q=p$ such that for all $i\in [q]$, $r\leq p_i-p_{i-1}\leq r+1$ and $f(p_{i-1}+1)< \dots <f(p_{i})$.  For $j\in [q]$, let $v_{p_{j-1}+1},\ldots,v_{p_j}$ be the \textit{$j$-th subsequence} of $P$.
\end{definition}

Given a properly ordered path $P=v_{p_0+1}\dots v_{p_1}v_{p_1+1}\cdots v_{p_2}\cdots v_{p_{q-1}+1}\cdots v_{p_q}$, we will say that the $j$-th subsequence, $v_{p_{j-1}+1}, \dotsc, v_{p_j}$, \textit{has type $z \in \mathbb{Z}^k$}
  if for $i \in [k]$, we have $z_i = 1$ when one of the vertices in the subsequence is in the part $V_i$
  and $z_i = 0$ otherwise. 
  From the definition of properly ordered, this means that $v_{p_{j-1}+1}, \dotsc, v_{p_j}$ has type $z \in \mathbb{Z}^k$
  if $z_i = |\{v_{p_{j-1}+1}, \dotsc, v_{p_j}\} \cap V_i|$ for
  every $i \in [k]$.

  It is clear that we need the parts which contain every $r$ consecutive vertices in $P$ to be distinct.  
  Given a properly ordered $(r-1)$-path, we will have this critical property if and only if the following condition is met for every $j\in [q-1]$, and $i\in\{p_{j-1} + 1, \ldots, p_{j}\}$, and $i'\in\{p_{j} + 1, \ldots, p_{j+1}\}$: 
  \begin{equation}\label{eq:unbalanced}
    \text{If $v_i$ and $v_{i'}$ are contained in the same part, then $i' - i \ge r$.}
  \end{equation}
  We can restate this observation in the following way:
  The parts which contain every $r$ consecutive vertices in $P$ are 
  distinct if and only if 
  for every $j\in [q-1]$ when we let $z$ be the type
  of the $j$-th subsequence and $z'$ be the type of the $(j+1)$-th subsequence
  we have the following: 
  \begin{equation}\label{eq:unbalanced_type}
    \text{For every $i \in [k]$, if $z_i = z'_i = 1$, then 
      $\sum_{v = i+1}^{k} z_v + \sum_{v = 1}^{i} z'_v$} \ge r.
  \end{equation}
  Note that if the $j$-th and $(j+1)$-th subsequences of $P$ both contain
  exactly $r$ vertices 
  (so, $\sum_{v = 1}^{k} z_v = \sum_{v = 1}^k z'_v = r$), then
  \eqref{eq:unbalanced_type} can be restated as the following:
  \begin{equation}\label{eq:unbalanced_typer}
    \text{For every $i \in [k]$, if $z_i = z'_i = 1$, then 
      $\sum_{v = 1}^{i} z_v 
      = r - \sum_{v = i+1}^{k} z_v 
      \le \sum_{v = i}^{i} z'_v$}.
  \end{equation}
  If the ordered pair $(z,z')$ satisfies \eqref{eq:unbalanced_type}, then
  we say that $(z,z')$ is \textit{valid}.

\begin{proof}[Proof of Lemma \ref{lem:seq}]

Let $2 \le r < k \le 2r - 1$ and let 
$\beta$ and $\sigma$ be constants such that
\begin{equation}\label{eq:constants}
  \frac{1}{n} \ll \beta \ll \sigma \ll \gamma \le \frac{1}{r}.
\end{equation}
Let $G$ be an $n$-vertex $k$-partite graph with ordered partition 
$\cP = (V_1, \dotsc, V_k)$ of $V = V(G)$ such that
\begin{equation}\label{eq:descending_part_sizes}
  \gamma n \le |V_k| \le |V_{k-1}| \le \dotsm \le |V_1| \le \frac{n}{r},
\end{equation}
and
\begin{equation}\label{eq:degree_condition}
  \delta_{\cP}(G) \ge 1 - \frac{1}{r} + \gamma.
\end{equation}

If $|V_1| \ge \frac{n}{r} - 2 \sigma n$, we define $1\le s \le k$ to be the largest integer such that 
$|V_s| \ge \frac{n}{r} - 2 \sigma n$; otherwise, we set $s=0$.

We start by greedily building a path $P_0'$ 
such that when $V' = V \setminus V(P_0')$ and
$V'_i = V_i \setminus V(P_0')$ for every $i \in [k]$, 
the following holds:
\begin{enumerate}[label=(T\arabic*),leftmargin=4\parindent]
  \item\label{T1} $|V'|$ is divisible by $r$;
  \item\label{T2} $|V'_i| = |V'|/r$ for every $i\in [s]$;
  \item\label{T3} $|V'_i| \ge \sigma n$ for every $i\in \{s + 1,\ldots,k\}$;
  \item\label{T4} $|V'_i| \le |V'|/r - \sigma n$ for every $i\in \{s+1,\ldots,k\}$; 
  \item\label{T5} $|V'| \ge (1 - 3r^2 \sigma)n$; and
  \item\label{T6} $P_0'$ is properly ordered and properly terminated.
\end{enumerate}

Let $z^{(0)}$ be the $(0,1)$-vector in $\mathbb{Z}^{k}$ 
in which the first $(r+1)$ entries are one and the remaining $k-r-1$ entries are zero.
For $j \in [r+1]$, let $z^{(j)}$ be $z^{(0)}$ minus the 
$j$-th standard basis vector, 
i.e., all of the last $k-r-1$ entries of $z^{(j)}$ are zero and all of the first $(r+1)$ entries of $z^{(j)}$ are one except for the $j$-th entry, which is zero.
Using \eqref{eq:unbalanced_type} and
\eqref{eq:unbalanced_typer} it is not hard to verify that
the following holds for every $j,j' \in [r+1]$:
\begin{enumerate}[label=(V\arabic*),leftmargin=4\parindent]
  \item\label{V1} $(z^{(0)}, z^{(j)})$ is valid;
  \item\label{V2} $(z^{(j)}, z^{(j')})$ is valid when $j \le j' + 1$; and
  \item\label{V3} $(z^{(j)}, z^{(j')})$ is not valid when $j \ge j' + 2$.
\end{enumerate}

Let $0 \le c_0 < r$ be such that $n - c_0$ is divisible by $r$ and
for $i \in [s]$, let 
\begin{equation}\label{eq:size_of_ci}
  c_i = \frac{n - c_0}{r} - |V_i|.
\end{equation}
Note, by \eqref{eq:descending_part_sizes} and the definition of $s$,
we have that $\frac{n}{r} - 2 \sigma n \le |V_s| \le \dotsm \le |V_1| \le \frac{n}{r}$, so
\begin{equation}\label{eq:trim1}
    2\sigma n \ge c_s \ge c_{s-1} \ge \cdots \ge c_1 \ge 0  .
\end{equation}

The sequences of vectors
\begin{equation*}
c_0 z^{(0)}, z^{(r+1)}, z^{(r)}, z^{(r-1)}, \dotsc, z^{(s+1)}, c_s z^{(s)}, 
c_{s-1} z^{(s-1)}, \dotsc, c_1 z^{(1)}, z^{(r+1)},
\end{equation*}
will serve as our template for $P_0'$.

That is, we greedily
build $P_0'$ so that 
\begin{itemize}
  \item the first $c_0$ subsequences are of type $z^{(0)}$ (these
    are the only subsequences that have $(r+1)$ instead of $r$ vertices);
  \item the next $(r-s+1)$ subsequences have types 
    $z^{(r+1)}, z^{(r)}, z^{(r-1)}, \dotsc, z^{(s+1)}$, respectively;
  \item the next $c_s$ subsequences are of type $z^{(s)}$,
    followed by $c_{s-1}$ subsequences of type $z^{(s-1)}$,
    \ldots, followed by $c_1$ subsequences of type $z^{(1)}$; and
  \item the last subsequence is of type $z^{(r+1)}$.
\end{itemize}
Note that it is possible to build $P_0'$ in this way by
\eqref{eq:degree_condition}, \eqref{eq:trim1}, 
\ref{V1} and \ref{V2}
(To see that \eqref{eq:trim1} is critical here, note that, by \ref{V3}, 
we need that if $t \in [s]$ is such that $c_t = 0$, then 
$c_{t-1} = c_{t-2} = \dotsm \ = c_1 = 0$.)
Define 
\begin{equation}\label{eq:size_of_q}
  q = c_0 + (r-s + 1) + \sum_{j = 1}^{s} c_j + 1
\end{equation}
and note that $q$ is the number of subsequences in $P_0'$.

\begin{claim}\label{clm:T1T6}
  The $P_0'$ constructed as described above satisfies conditions~\ref{T1}--\ref{T6}.
\end{claim}

\begin{proof}[Proof of Claim~\ref{clm:T1T6}]~\\
\noindent\textbf{\ref{T6}:}
The construction of $P_0'$ requires $P_0'$ to be
properly ordered and properly terminated (even when $c_0 = 0$).

\noindent\textbf{\ref{T1}:}
Recall that each subsequence has $r$ vertices except the first $c_0$, which have $r+1$. By \eqref{eq:size_of_q}, the number of vertices
in $P_0'$ is 
\begin{equation}\label{eq:size_of_P0}
  p = c_0(r+1) + (r-s+1)r + \sum_{j=1}^{s} c_j r + r
   = c_0 + qr.
\end{equation}
So, since $n - c_0$ is divisible by $r$, we have that
$|V'|=n-p$ is divisible by $r$.

\noindent\textbf{\ref{T5}:}
By \eqref{eq:trim1}, \eqref{eq:size_of_q}, \eqref{eq:size_of_P0} 
and the fact that $s\leq r$ and $c_0 < r$, we have that 
\begin{equation}\label{eq:size_of_p}
  p = c_0 + qr = c_0(r+1) + (r-s+2)r + r\sum_{j=1}^s c_j \le 3 \sigma r^2 n . 
\end{equation}

\noindent\textbf{\ref{T3}:}
By \eqref{eq:descending_part_sizes}, for all $i\in\{s+1,\ldots,k\}$,
\begin{equation*}
  |V_i'| = |V_i\setminus V(P_0')| \ge \gamma' n - 3\sigma r^2n \ge \sigma n . 
\end{equation*}

\noindent\textbf{\ref{T2}:}
By \eqref{eq:size_of_ci} and \eqref{eq:size_of_P0}, for all $i \in [s]$, 
\begin{equation*}
  |V'_i| = |V_i| - q + c_i 
  = |V_i| - q + \frac{n - c_0}{r} - |V_i|
  = \frac{n - c_0}{r} - \frac{p - c_0}{r} 
  = \frac{n - p}{r}
  = \frac{|V'|}{r} .
\end{equation*}

\noindent\textbf{\ref{T4}:}
Consider two cases:
If $s+1 \le i \le r$, then, because every subsequence of $P_0$
except exactly one intersects $V_i$, we have
\begin{equation*}
  |V'_i| = |V_i| - q + 1 =
  |V_i| - \frac{p - c_0}{r} + 1 < 
  \left(\frac{n}{r} - 2 \sigma n\right) - \frac{p}{r} + 2 
  = \frac{|V'|}{r} - 2 \sigma n + 2 < \frac{|V'|}{r} - \sigma n .
\end{equation*}
If $r+1 \le i \le k$, \eqref{eq:descending_part_sizes} implies
that $|V_i| \le n/i \le n/(r+1)$, so with \eqref{eq:size_of_p},
\begin{equation*}
  |V'_i| \le |V_i| \le \frac{n}{r+1} = \frac{n}{r} - \frac{n}{r(r+1)} \le 
\frac{n}{r} - 3\sigma r n  - \sigma n \le \frac{n}{r} - \frac{p}{r} - \sigma n =
\frac{|V'|}{r}  - \sigma n.
\end{equation*}

\noindent This concludes the proof of Claim~\ref{clm:T1T6}.
\end{proof}

Now we consider \ref{A1}.  We stress that the issue here is largely numerical, which explains the general nature of the next two claims.  Claim \ref{clm:S2C1} provides the template and Claim \ref{clm:S2C2} shows that $V'$ can be partitioned according to the template so that \ref{A1} holds.  The purpose of partitioning according to this specific template is to set things up so that \ref{A3} and \ref{A4} will be able to be satisfied in the end.

Let $\cZ$ be the set of $(0,1)$-vectors in $\mathbb{Z}^k$ 
such that the first $s$ entries are one and exactly $r-s$ of the 
remaining $k-s$ entries are one (so, for every $z \in \cZ$ 
exactly $r$ of the $k$ entries of $z$ are one and the remaining
$k-r$ entries are zero).
Note that $\ell = \binom{k-s}{r-s}$ is the order of $\cZ$.

\begin{claim}\label{clm:S2C1}
  There exists a $k \times \ell$ $(0,1)$-matrix $A = [a_{i,j}]$ such
  that the $\ell$ columns of $A$ are the vectors in $\cZ$ where
  the columns of $A$ are ordered so that 
  \begin{itemize}
    \item the first column is $(\underbrace{1,\dotsc,1}_\text{$r$ times}, \underbrace{0,\dotsc,0}_\text{$k-r$ times})^T$;
    \item the last column is $(\underbrace{1,\dotsc,1}_\text{$s$ times}, \underbrace{0,\dotsc,0}_\text{$k-r$ times}, \underbrace{1,\dotsc,1}_\text{$r-s$ times})^T$; and 
    \item for every $j \in [\ell - 1]$ and $i \in [k]$, 
      \begin{equation}\label{eq:S2P1}
        \text{if $a_{i,j} = a_{i, j+1} = 1$, then $\sum_{v = 1}^{i} a_{v,j} \le \sum_{v=1}^{i} a_{v,j+1}$} \qquad \text{(c.f.\ \eqref{eq:unbalanced_typer}).}
      \end{equation}
  \end{itemize}
\end{claim}

\begin{proof}[Proof of Claim~\ref{clm:S2C1}]
  The proof is by induction on $k-s$. 
  Note that if either $k=r$ or $r=s$, then the claim is trivially true.
  In particular, this establishes the base case since $k-s = 0$ implies $k=r=s$. 
  Now suppose that $k > r > s$.
  Let $\cZ'$ be the vectors in $\cZ$ in which the $(s+1)$-th entry
  is one and let $\cZ'' = \cZ \setminus \cZ'$.
  Let $\ell' = |\cZ'| = \binom{k-s-1}{r-s-1}$
  and $\ell'' = |\cZ''| = \binom{k-s-1}{r-s}$.
  By the induction hypothesis (with $k$, $r$, and $s+1$ playing the roles of $k$, $r$, and $s$, respectively), 
  we can populate the first $\ell'$ columns of $A$ with the vectors in $\cZ'$
  so that the first column is $(\underbrace{1,\dotsc,1}_\text{$r$ times}, \underbrace{0,\dotsc,0}_\text{$k-r$ times})^T$,
  the $\ell'$-th column is $(\underbrace{1,\dotsc,1}_\text{$s+1$ times}, \underbrace{0,\dotsc,0}_\text{$k-r$ times}, \underbrace{1,\dotsc,1}_\text{$r-s-1$ times})^T$, and
  \eqref{eq:S2P1} holds for $j\in[\ell' - 1$].
  Similarly, by the induction hypothesis (with $k - s - 1$, $r-s$, and $0$ playing the roles of $k$, $r$, and $s$, respectively),
  we can populate the remaining columns of $A$ with $\cZ''$
  so that the $(\ell'+1)$-th column is $(\underbrace{1,\dotsc,1}_\text{$s$ times}, 0, \underbrace{1,\dotsc,1}_\text{$r-s$ times}, \underbrace{0,\dotsc,0}_\text{$k-r-1$ times})^T$,
  the last column is $(\underbrace{1,\dotsc,1}_\text{$s$ times}, \underbrace{0,\dotsc,0}_\text{$k-r$ times}, \underbrace{1,\dotsc,1}_\text{$r-s$ times})^T$, 
  and \eqref{eq:S2P1} holds for $\ell' + 1 \le j \le \ell - 1$.
  The claim then follows because \eqref{eq:S2P1} holds when $j=\ell'$.
\end{proof}

Let $A$ be the matrix guaranteed by Claim~\ref{clm:S2C1}.
\begin{claim}\label{clm:S2C2}
  Let $b = (|V'_1|, |V'_2|, \dotsc, |V'_k|)^T$.
  There exists $x \in \mathbb{Z}^{\ell}$ such that 
  $x_j \ge \beta n$ for every $j \in [\ell]$ and such that $A x = b$.
\end{claim}

\begin{proof}[Proof of Claim~\ref{clm:S2C2}]
  We will iteratively construct a sequence of vectors $x^{(0)}, x^{(1)}, \dotsc, x^{(T)} \in \mathbb{Z}^\ell$ such that
  $x = x^{(T)}$ meets the conditions of the claim.
  For $t \ge 0$, define $b^{(t)} = b - A x^{(t)}$; $n^{(t)} = \sum_{i=1}^{k}b^{(t)}_i$;
  and the following properties:
  \begin{enumerate}[label=(P\arabic*),leftmargin=4\parindent]
    \item\label{P1} $n^{(t)} \ge 0$ is divisible by $r$;
    \item\label{P2} $b^{(t)}_i = n^{(t)}/r$ for every $1 \le i \le s$;
    \item\label{P3} $0 \le b^{(t)}_i \le n^{(t)}/r$ for every $s+1 \le i \le k$; and
    \item\label{P4} $x^{(t)}_j \ge \beta n$ for every $j \in [\ell]$.
  \end{enumerate}

  To begin the construction, we let 
  $m = \ceiling{\beta n}$ and $x^{(0)}_j = m$ for every $j \in [\ell]$. 
  Clearly, we have that \ref{P4} holds for $t = 0$.
  First note that $n^{(0)} = |V'| - r \ell m$ so, by \ref{T1}, we have that \ref{P1} holds for $t = 0$.
  By \ref{T2}, we also have that 
  \begin{equation*}
    b^{(0)}_i = |V'|/r - \ell m = n^{(0)}/r \qquad \text{for every $i \in [s]$},
  \end{equation*}
  so \ref{P2} holds for $t = 0$.
  By \ref{T3}, we have $b^{(0)}_i \ge b_i - \ell m \ge \sigma n - \ell m > 0$ for every $s+1 \le i \le k$, 
  and with \ref{T4} we have that 
  \begin{equation*}
    b^{(0)}_i \le |V'|/r - \sigma n \le  |V'|/r - \ell m = n^{(0)}/r \qquad \text{for every $s+1 \le i \le k$}.
  \end{equation*}
  Therefore, \ref{P3} also holds for $t = 0$.

  Now assume \ref{P1}, \ref{P2}, \ref{P3}, and \ref{P4} hold for some $t \ge 0$.
  If $b^{(t)}_i = 0$ for every $i \in [k]$, then $A x^{(t)} = b$, so with \ref{P4}, we 
  can let $t = T$ and end the construction, because $x = x^{(t)} = x^{(T)}$ meets the conditions of the claim. 
  Otherwise, let $I = \{i \in [k]: b^{(t)}_i = n^{(t)}/r\}$.
  Note that \ref{P2} implies that $[s] \subseteq I$ and by 
  \ref{P3} we have that $b^{(t)}_i \le n^{(t)}/r - 1$ for every $i \in [k] \setminus I$.
  We clearly have that $|I| \le r$ and, by \ref{P1}, \ref{P2} and \ref{P3}, 
  there exists $I \subseteq I' \subseteq [k]$ such that $|I'| = r$ and $b^{(t)}_i > 0$ for every $i \in I'$.
  Now let $j^{(t)}$ be the column of $A$ such that $a_{i,j^{(t)}} = 1$ if and only if $i \in I'$.
  If we then let 
  \begin{equation*}
    x^{(t+1)}_j = 
    \begin{cases} 
      x^{(t)}_j + 1 & \text{if $j=j^{(t)}$} \\
      x^{(t)}_j     & \text{otherwise}
    \end{cases}
  \end{equation*}
  it is clear that \ref{P1}, \ref{P2}, \ref{P3}, and \ref{P4} all hold with $t$ set to $t+1$.
\end{proof}

Now we use the preceding claims to show that \ref{A1} and \ref{A2} hold.  Let $i \in [k]$ and recall that, since $Ax = b$, we have 
$\sum_{j=1}^{\ell} a_{i,j} \cdot x_j = b_i = |V'_i|$. 
Therefore, for every $i \in [k]$, 
we can uniformly at random select a partition of $V'_i$ into
$\ell$ parts $V(i,1), \dotsc, V(i,\ell)$ so that for every
$j \in [\ell]$, we have $|V(i,j)| = a_{i,j} \cdot x_j$. 
(Note that we are allowing parts to be empty in these partitions).

Let $j \in [\ell]$, 
and 
note that since exactly $r$ entries in the $j$-th column of $A$ are $1$, 
there exists a unique sequence $1 \le i_{j,1} < \dotsm < i_{j,r} \le k$ such that
$|V(i_{j,1},j)| = \dotsm = |V(i_{j,r},j)| = x_j \ge \beta n$,
so \ref{A1} holds.

Let $\cP_j = (V(i_{j,1},j), \dotsc, V(i_{j,r},j))$ and
$G_j = G[ V(i_{j,1}, j) \cup \dotsm \cup V(i_{j,r},j)]$.
Therefore, \eqref{eq:degree_condition}, \ref{T5}, and the 
Chernoff and union bounds imply that, with high probability,
there exists an outcome where for every 
$j \in [\ell]$, $h \in [r]$, and $v \in V \setminus V_{i_{h,j}}$ we have that
\begin{equation}\label{eq:into_Vij}
  \deg(v, V(i_{j,h},j)) \ge \left(1 - \frac{1}{r} + \frac{\gamma}{2}\right)|V(i_{j,h},j)|.
\end{equation}
Fix such an outcome.
Note that $V(P'_0), V(G_1), \dotsc, V(G_\ell)$ is a partition of $V(G)$,
and for every $j \in [\ell]$, $G_j$ is a balanced $r$-partite graph with
ordered partition $\cP_j$ such 
that each part has order at least $\beta n$ and,
with \eqref{eq:into_Vij} we have
\begin{equation}\label{eq:G_j_partial_degree}
  \delta_{\cP_j}(G_j) \ge 1 - \frac{1}{r} + \frac{\gamma}{2},
\end{equation}
i.e., \ref{P2} holds.

To see that \ref{A4} holds,  
note that, by 
the ordering of the columns
of $A$ (c.f. \eqref{eq:S2P1}) and \eqref{eq:into_Vij}, 
we can greedily construct $\ell-1$ 
vertex disjoint $(r-1)$-paths $P_1, \dots, P_{\ell-1}$
each on exactly $2r$ vertices so that, for every $j \in [\ell-1]$, 
the initial $r$ vertices of $P_j$ respect the sequence $\mathcal{P}_j$
and the final $r$ vertices of $P_j$ respect the sequence $\mathcal{P}_{j+1}$.

Finally, we now show that \ref{A3} holds.  
To see this, first note that, by Claim~\ref{clm:S2C1}, 
if we let $z$ be the last column of $A$ and 
$z'$ be the first column of $A$, then by \eqref{eq:unbalanced_typer},
we have that $(z,z')$ is valid.
Therefore, since \ref{T6} implies that the $(r-1)$-path $P_0'$ is 
properly terminated,
we can use \eqref{eq:into_Vij} to greedily prepend
$r$ vertices to $P_0'$ to create an $(r-1)$-path in which 
the initial $r$ vertices respect the sequence $\mathcal{P}_\ell$
while avoiding the path $P_{\ell-1}$.
Again because $P_0'$ is properly terminated, \eqref{eq:into_Vij} implies
that we can greedily append $r$ vertices to this path to create an 
$(r-1)$-path $P_0$
so that the final $r$ vertices
of $P_0$ respect the sequence $\mathcal{P}_1$
and so that $P_0$ avoids $P_1$.
This completes the proof.
\end{proof}

\section{Conclusion}\label{sec:conclusion}

\subsection{Exact version}

The main open problem which remains is to prove an exact version of Theorem \ref{thm:main}.  Note that it is possible that in the unbalanced case, there are extra variants of Catlin's example.

\subsection{Total degree version}

Another direction is to consider minimum total degree conditions for perfect $K_r$-tilings and Hamiltonian $(r-1)$-paths.  In this direction, Johansson, Johansson, and Markstr\"{o}m \cite{JJM} proved that if $G$ is a balanced 3-partite graph on $n$ vertices with $\delta(G)\geq n/2$, then $G$ has a perfect $K_3$-tiling.  Later, Lo and Sanhueza-Matamala \cite{LS} proved that if $G$ is a balanced $r$-partite graph on $n$ vertices with $\delta(G)\geq \left(1-\frac{3}{2r}+o(1)\right)n$, then $G$ has a perfect $K_r$-tiling, which is asymptotically best possible.  

It would be interesting to study the unbalanced version of this result and extend it to Hamiltonian $(r-1)$-cycles.  This was done for $r=2$ in \cite{DKPT}, but the degree condition is quite complicated (in some sense necessarily so, since it is asymptotically tight in all cases) and thus determining an asymptotically tight minimum degree condition for perfect $K_r$-tilings in all valid $k$-partite graphs seems challenging.

As a start, we conjecture the following sufficient condition for perfect $K_r$-tilings (which will be asymptotically necessary in certain cases).

\begin{conjecture}
Let $k\geq r\geq 2$ and $\gamma>0$.  If $G$ is a $k$-partite graph with all parts at most $n/r$ and $\delta(V_i)\geq \left(1-\frac{1}{2r}+\gamma\right)n-|V_i|$ for all $i\in [k]$, then $G$ has a perfect $K_r$-tiling.  
\end{conjecture}

{\bf Acknowledgement.}
 We  thank the anonymous referee for the helpful comments which improved the presentation of this paper.

\bibliographystyle{alpha}

\end{document}